\documentclass[notitlepage,12pt,draft]{article}
\usepackage{amssymb,amsfonts,amsmath,geometry,esint,comment,upgreek,cases,stmaryrd,ulem}

\catcode`\@=11
\@addtoreset{equation}{section}

\catcode`\@=12

\usepackage{color}

\usepackage{scalerel,stackengine}
\stackMath
\newcommand\reallywidehat[1]{%
\savestack{\tmpbox}{\stretchto{%
  \scaleto{%
    \scalerel*[\widthof{\ensuremath{#1}}]{\kern-.6pt\bigwedge\kern-.6pt}%
    {\rule[-\textheight/2]{1ex}{\textheight}}
  }{\textheight}%
}{0.5ex}}%
\stackon[1pt]{#1}{\tmpbox}%
}
\def\R{\mathbb{R}}

\def\f{\varphi}

\def\div{{\rm div}}

\def\me{\mathrm{\bf e}}
\def\is{{\lfloor{s}\rfloor}}
\def\istwo{{\lfloor{2s}\rfloor}}
\def\isalpha{{\lfloor{\alpha}\rfloor}}

\def\intr{\int\limits_{-\infty}^{+\infty}}

\def\Ds{\left(-\Delta\right)^{\!s}} 

\def\Db{{\mathcal D}_{b}}

\def\Db{\mathbb D}
\def\Lb{\mathbb L}

\newcommand{\bL}{\mathbb L}

\def\proof{\noindent{\textbf{Proof. }}}
\def\QED{\hfill {$\square$}\goodbreak \medskip}

\newtheorem{Theorem}{Theorem}[section]
\newtheorem{Lemma}[Theorem]{Lemma}
\newtheorem{Proposition}[Theorem]{Proposition}
\newtheorem{Corollary}[Theorem]{Corollary}
\newtheorem{Remark}[Theorem]{Remark}

\usepackage{mathtools}

\begin{document}

\title 
{Fractional operators as traces \\of  operator-valued curves}

\author{Roberta Musina\footnote{Dipartimento di Scenze Matematiche, Informatiche e Fisiche, Universit\`a di Udine, via delle Scienze, 206 -- 33100 Udine, Italy. Email: {roberta.musina@uniud.it}.}~ 
and \setcounter{footnote}{6}
Alexander I. Nazarov\footnote{
St.Petersburg Dept. of Steklov Institute, Fontanka 27, St.Petersburg, 191023, Russia, 
and St.Petersburg State University, 
Universitetskii pr. 28, St.Petersburg, 198504, Russia. E-mail: al.il.nazarov@gmail.com. 
{Partially supported by RFBR grant 20-01-00630.}
}
}
\date{}

\maketitle

\noindent
{\small {\bf {Abstract.}}  
We  relate non integer powers ${\mathcal L}^{s}$, $s>0$ of a given (unbounded) positive self-adjoint operator $\mathcal L$ in a real separable Hilbert space $\mathcal H$
with a certain differential operator of  order $2\lceil{s}\rceil$, acting on even curves $\R\to \mathcal H$. This extends the results by Caffarelli--Silvestre and Stinga--Torrea regarding the characterization
of fractional powers of differential operators via an extension problem.
}

\medskip
\noindent
{\footnotesize {\bf Keywords:}  Higher order fractional operators; Operator-valued functions.}

\noindent
{\footnotesize {\bf 2010 Mathematics Subject Classification:} 34G10; 46E40; 35R11.}  

\normalsize

\section{Introduction} 

The study of the fractional powers of differential operators via their relations with generalized harmonic extensions
and corresponding Dirichlet-to-Neumann operators began more than fifty years ago \cite{MO} and became popular thanks to the celebrated work \cite{CS} of Caffarelli and Silvestre, which stimulated a fruitful line of research. The idea of relating the operators  $\Ds $, $s\in(0,1)$, acting on $\R^n$ and $-\div(y^{1-2s}\nabla)$ acting on $\R^n\times\R_+$,
has been adapted to cover much more general situations. The first contribution in this direction is due to Stinga and Torrea \cite{ST};  important generalizations were given in  \cite{AEW, GMS}. 

The case of higher order powers of $\Ds$ has been investigated firstly in \cite{CdMG} via conformal geometry techniques.
We also cite  \cite{CC, MMG, GR, JX, RS}, the more recent papers \cite{Case, CM} and references there-in. 

Before describing our results, let us notice that any extension $w=w(\cdot,y)$ of a given $u=u(\cdot)$ can be related with 
the curve $y\mapsto w(\cdot, y)$ taking values in a suitable  function space. 
In the present paper we use this interpretation to handle any non-integer power $s>0$ of a linear operator $\mathcal L$ in quite a  general framework.

Let ${\mathcal H}$ be a separable real Hilbert space with scalar product $(\cdot,\cdot)_{\!\mathcal H}~\!$ and norm $\|\cdot\|_{\!\mathcal H}~\!$. Let
$$
{\mathcal L}: \mathcal D({\mathcal L})\to {\mathcal H} ~,\qquad  \mathcal D({\mathcal L})\subseteq {\mathcal H}
$$
be a given unbounded,  self-adjoint operator. In order to simplify the exposition, we first assume that 
$\mathcal L$ is positive definite and has discrete spectrum  (some generalizations are given in Section  \ref{S:generalization}). 
We organize the spectrum  of ${\mathcal L}$ 
in a nondecreasing sequence of eigenvalues $(\lambda_j)_{j\ge 1}$, counting with their multiplicities, and denote by $\f_j\in \mathcal D({\mathcal L})$ a complete orthonormal system of corresponding eigenvectors.

Given $s\in\R$, the $s$-th power of $\mathcal L$ in the sense of spectral theory
is the operator
\begin{equation}
\label{eq:dual0}
{\mathcal L}^{s} u= \sum_{j=1}^\infty \lambda_j^{s}u_j~\!\f_j~\!,\quad\text{where}\quad u_j=(u,\f_j)_{\!\mathcal H}~\!,
\end{equation}
so that $\mathcal L^{0}$ is the identity in $\mathcal H$.
If $s>0$, the natural domain of the quadratic form 
$$
u\mapsto ({\mathcal L}^{s}u,u)_{\!\mathcal H}=\sum_{j=1}^\infty\lambda_j^{s}u_j^2
$$
is denoted by
$\mathcal H^s_{\mathcal L}$. Clearly $\mathcal H^s_{\mathcal L}$ coincides with the domain of ${\mathcal L}^{\frac{s}{2}}$; it is a Hilbert space with scalar product and norm
given by
\begin{equation}
\label{eq:Hs_norm}
(u,v)_{\mathcal H^s_{\mathcal L}}= ({\mathcal L}^{\frac{s}{2}}u,{\mathcal L}^{\frac{s}{2}}v)_{\!\mathcal H}~,\quad  \|u\|_{\mathcal H^s_{\mathcal L}} =\|{\mathcal L}^{\frac{s}{2}}u\|_{\!\mathcal H}~\!.
\end{equation}
We identify the dual space $(\mathcal H^{s}_{\mathcal L})'$ with $\mathcal H^{-s}_\mathcal L=\{\mathcal L^s u~|~u\in \mathcal H^s_{\mathcal L}~\}$ 
via the identity
$$
\langle {\mathcal L}^s u,v\rangle=({\mathcal L}^\frac{s}{2}u, {\mathcal L}^\frac{s}{2}v)_{\!\mathcal H}~\quad \text{for any $u,v\in \mathcal H^s_{\mathcal L}$.}
$$
Notice that $\mathcal L^{s}$ is an isometry $\mathcal H^{s}_\mathcal L\to \mathcal H^{-s}_\mathcal L$
with inverse $\mathcal L^{-s}$.

\medskip
In this paper we relate the operator ${\mathcal L}^{s}: \mathcal H^s_{\mathcal L}\to \mathcal H^{-s}_{\mathcal L}$ for $s>0$ non-integer to certain linear
operator acting on {\bf even} curves $\R\to \mathcal H^s_{\mathcal L}$ 
(this simplifies the treatment in case of higher powers $s>1$, compare with  \cite{CM}). 

\medskip

Let $b\in(-1,1)$.
Denote by $L^{2;b}(\R\to\mathcal H)$  
the Hilbert space of curves $U:\R\to\mathcal H$ such that $\|U(y)\|^2_{\!\mathcal H}$ is integrable on $\R$ with respect to the measure $|y|^bdy$. Further,
$L^{2;b}_\me(\R\to\mathcal H)$ stands for the subspace of even curves.

For $U\in L^{2;b}_\me(\R\to\mathcal H)$ we consider the (unbounded) operators
\begin{equation}
\label{eq:Db}
\Db_bU=-\partial^2_{yy}U-by^{-1}\partial_yU=-|y|^{-b}\partial_y(|y|^b\partial_yU)~,\qquad \Lb_bU=\Db_bU+{\mathcal L}U.
\end{equation} 

Denoting by $U_j(y)=(U_j(y),\f)_{\!\mathcal H}$ the coordinates of $U(y)$, we have
$$
\Lb_bU=\sum_{j=1}^\infty\big((\Db_b+\lambda_j)U_j\big)~\!\f_j,
$$
and the corresponding quadratic form reads
$$
(\Lb_bU,U)_{L^{2,b}}\!\!=\!\!\intr\!\!|y|^b\big(\|\partial_y U(y)\|_{\mathcal H}^2+\|\mathcal L^\frac12(U(y))\|_{\mathcal H}^2)dy=
\sum_{j=1}^\infty\intr\!\!|y|^b(|\partial_yU_j|^2+\lambda_j|U_j|^2)dy.
$$

In Section \ref{SS:U_space} we study in detail the 
natural domain
$$
H^{k;b}_{\mathcal L,\me}(\R\to\mathcal H)\subset L^{2;b}_\me(\R\to\mathcal H)~,\qquad k\in\mathbb N,
$$
of the quadratic form 
$U\mapsto (\Lb_b^{k}U,U)_{L^{2;b}}$.  Lemma \ref{def:Lspace}
 provides explicit expressions for its Hilbertian scalar product and related norm, which are denoted by 
$(\cdot, \cdot)_{\!H^{k;b}_{\mathcal L,\me}}$, $\|\cdot\|_{\!H^{k;b}_{\mathcal L,\me}}^2$,
respectively, and shows that the Dirac-type trace function $\delta_0(V)=V(0)$ is continuous 
from $H^{k;b}_{\mathcal L,\me}(\R\to\mathcal H)$ into $\mathcal H^{k-\frac{1+b}{2}}_\mathcal L$.

Our main results involve the linear transform 
\begin{equation}
\label{eq:best}
\mathcal P_{\!s}[u](y)= \frac{2^{1-s}}{\Gamma(s)} 
\sum_{j=1}^\infty ({\sqrt{\lambda_j}}|y|)^s K_s(\sqrt{\lambda_j}|y|)~\!  u_j\f_j
\end{equation}
for $u=\sum_j u_j\f_j\in \mathcal H$  and $y\in \R$, where $K_s$ is the modified Bessel function of the second kind
(the  Macdonald function; compare with  \cite{ST}, where $s\in(0,1)$ is assumed).

Due to the regularity and decaying properties of the Bessel functions, in Lemma \ref{L:P_smooth} of Appendix \ref{A}, we prove that for any $u\in \mathcal H$,
$\mathcal P_{\!s}[u]$ is an even curve in $\mathcal H$; in addition
$\mathcal P_{\!s}[u]\in {\mathcal C}^\infty(\R_+\to \mathcal H^\sigma_\mathcal L)$ for every $\sigma>0$.

To state our main result we introduce the floor and ceiling notation
$$
\lfloor{s}\rfloor\   \mbox{$:=$ integer part of $s$};\quad \lceil{s}\rceil:=\lfloor{s}\rfloor+1.
$$

\begin{Theorem}
\label{T:draft} 
Let $s>0$ be non-integer. We put
$$
{\mathfrak b}:=1-2(s-\lfloor{s}\rfloor)\in(-1,1)~\!.
$$
For any $u\in \mathcal H^s_{\mathcal L}$ the following facts hold.
\begin{enumerate}
\item[$i)$] 
\begin{equation}
\label{eq:isometry}
\|\mathcal P_{\!s}[u]\|_{\!H^{\lceil{s}\rceil;{\mathfrak b}}_{\mathcal L,\me}}^2=2{d_s}\|u\|_{\mathcal H^s_{\mathcal L}}^2\quad\text{where}
\quad d_s=2^{\mathfrak 
b}\Gamma\Big(\frac{1+{\mathfrak 
b}}{2}\Big)\frac{\lfloor{s}\rfloor!}{\Gamma(s)}~. 
\end{equation}
That is, up to a  constant, the transform
$\mathcal P_{\!s}$ is an isometry $\mathcal H^s_{\mathcal L}\to 
H^{\lceil{s}\rceil;{\mathfrak b}}_{\mathcal L,\me}(\R\to\mathcal H)$;
\item[$ii)$]~$\mathcal P_{\!s}[u]$ achieves
\begin{equation}
\label{eq:U_minimum}
\min_{U\in H^{\lceil{s}\rceil;{\mathfrak b}}_{\mathcal L,\me}(\R\to\mathcal 
H)\atop U(0)=u} \|U\|_{\!H^{\lceil{s}\rceil;{\mathfrak b}}_{\mathcal 
L,\me}}^2=2d_s\|u\|^2_{\mathcal H^s_\mathcal L}; 
\end{equation}
\item[$iii)$]~$(\mathcal P_{\!s}[u],V)_{\!H^{\lceil{s}\rceil;{\mathfrak 
b}}_{\mathcal L,\me}}=~\!2d_s~\!\langle{\mathcal L}^{s}u,V(0)\rangle$ for any 
$V\in H^{\lceil{s}\rceil;{\mathfrak b}}_{\mathcal L,\me}(\R\to\mathcal H)$;
\item[$iv)$]~$\mathcal P_{\!s}[u]$ solves the differential equation
\begin{equation}
\label{eq:ODE}
\Lb_{\mathfrak b}^{\lceil{s}\rceil}\mathcal P_{\!s}[u]=0~~ \text{in~ $\R_+$}
\end{equation}
and satisfies
\begin{equation}
\label{eq:bdry1}
\lim\limits_{y \to 0^+}\mathcal P_{\!s}[u](0)=
u\quad\text{in $\mathcal H^{s}_{\mathcal L}$}~,\quad \lim\limits_{y \to 
0^+}y^{\mathfrak b}~\!{\partial_y} \big(\bL_{\mathfrak 
b}^{\lfloor{s}\rfloor}\mathcal P_{\!s}[u]\big)(y)= - d_s\,{\mathcal L}^{s} u
\quad\text{in $\mathcal H^{-s}_{\mathcal L}$.}
\end{equation}
\end{enumerate}
\end{Theorem}
Additional information on the regularity  of $\mathcal P_{\!s}[u]$ and on its behavior at $\{y=0\}$ is given in  Appendix \ref{A},
see in particular Theorems \ref{C:3} and  \ref{T:4}.  Corollary \ref{C:convergence}  improves the convergence in
\cite[Theorem 1.1]{ST}, where $s\in(0,1)$ is assumed; in Subsection \ref{SS:isometric} we point out some isometric 
properties of the operator $P_{\!s}$ in the spirit of \cite{2015}.

We can also consider negative, non-integer orders. 

Let $s>0$. 
If $\zeta\in \mathcal H^{-s}_{\mathcal L}$ then $\mathcal L^{-s}\zeta\in  \mathcal H^{s}_{\mathcal L}$, so that for any
$y\in \R$ we can compute
$$
\mathcal P_{\!s}[\mathcal L^{-s}\zeta](y)=
\frac{2^{1-s}}{\Gamma(s)} \sum_{j=1}^\infty \lambda_j^{-s}\zeta_j~\! ({\sqrt{\lambda_j}}|y|)^s 
K_s(\sqrt{\lambda_j}|y|)\f_j.
$$
The next result is in fact a corollary of Theorem \ref{T:draft}.

\begin{Theorem}
\label{T:negative} 
Let $s>0$, ${\mathfrak b}\in(-1,1)$ be as in Theorem \ref{T:draft}. For any 
$\zeta\in \mathcal H^{-s}_{\mathcal L}$ the 
following facts hold.
\begin{enumerate}
\item[$i)$]~
\begin{equation}
\label{eq:isometry_negative}
\|\mathcal P_{\!-s}[\zeta]\|_{\!H^{\lceil{s}\rceil;{\mathfrak b}}_{\mathcal 
L,\me}}^2=2{d_s}\|\zeta\|_{\mathcal H^{-s}_{\mathcal L}}^2
\quad\text{where} \quad \mathcal P_{\!-s}:=\big(\mathcal P_{\!s}\circ\mathcal L^{-s}\big)~\!.
\end{equation}
That is, up to a constant, the transform $\mathcal P_{\!-s}$
is an isometry $\mathcal H^{-s}_{\mathcal L}\to 
H^{\lceil{s}\rceil;{\mathfrak b}}_{\mathcal L,\me}(\R\to\mathcal H)$;
\item[$ii)$]~$\mathcal P_{\!-s}[\zeta]$  achieves
\begin{equation}
\label{eq:U_minimum_negative}
\min_{U\in H^{\lceil{s}\rceil;{\mathfrak b}}_{\mathcal L,\me}(\R\to\mathcal 
H)} \big(\|U\|_{\!H^{\lceil{s}\rceil;{\mathfrak b}}_{\mathcal L,\me}}^2-
4d_s\langle\zeta, U(0)\rangle\big)=-2d_s\|\zeta\|_{\mathcal H^{-s}_{\mathcal L}}^2.
\end{equation}
\item[$iii)$]~$(\mathcal P_{\!-s}[\zeta],V)_{\!H^{\lceil{s}\rceil;{\mathfrak 
b}}_{\mathcal L,\me}}
=~\!2d_s~\!\langle\zeta,V(0)\rangle$ for any $V\in 
H^{\lceil{s}\rceil;{\mathfrak b}}_{\mathcal L,\me}(\R\to\mathcal H)$;
\item[$iv)$]~$\mathcal P_{\!-s}[\zeta]$ solves  the differential equation
$$
\Lb_{\mathfrak b}^{\lceil{s}\rceil}\mathcal P_{\!-s}[\zeta]=0~~ \text{in~ 
$\R_+$}
$$
and satisfies
$$
\lim\limits_{y \to 0^+}y^{\mathfrak b}~\!{\partial_y} 
\big(\bL_{\mathfrak b}^{\lfloor{s}\rfloor}\mathcal P_{\!-s}[\zeta]\big)(y)= - 
d_s\,\zeta~~ \text{in~ $\mathcal H^{-s}_{\mathcal L}$}~,\quad
\lim\limits_{y\to 0^+}\mathcal P_{\!-s}[\zeta](y)=\mathcal L^{-s}\zeta~~ \text{in~ $\mathcal H^{s}_{\mathcal L}$}~\!.
$$
\end{enumerate}
\end{Theorem}

\medskip

{The paper is organized as follows.}
We start by introducing and studying, in Section \ref{S:spaces}, some Sobolev-type spaces $H^{k;b}_\me(\R)$  depending on 
the integer $k\ge 1$ and on the parameter $b\in(-1,1)$. In Section \ref{S:Bessel} we investigate the properties of the functions
\begin{equation}
\label{eq:psi_notation}
\psi_s(y)=c_s |y|^s K_s(|y|), \quad  c_s= \frac{2^{1-s}}{\Gamma(s)},
\end{equation}
which are involved in the definition of the operator $u\mapsto \mathcal P_{\!s}[u]$. The main result here is Theorem \ref{T:Aexistence}, which constitutes the basic tool in the proof of Theorem \ref{T:draft}. 

Section \ref{SS:U_space} contains the description of the Hilbert space $H^{k;b}_{\mathcal L,\me}(\R\to\mathcal H)$
of even curves in $\mathcal H$ mentioned above, and  the proofs of Theorems \ref{T:draft} and \ref{T:negative}.

Generalizations and examples are given in 
Section \ref{S:generalization}. 

As already mentioned, the Appendix contains several results about the operator $\mathcal P_{\!s}$.

\bigskip
\noindent
{\bf Acknowledgment.}
We wish to thank Prof. Fedor Petrov for having suggested us the proof of Lemma \ref{L:thanks} and Prof.
Nicolay Filonov for the useful and stimulating discussions.\medskip

\footnotesize{
\noindent
{\bf Notation.} 
Let $X$ be a Hilbert space with scalar product $(\cdot,\cdot)_X$ and norm $\|\cdot\|_X$. For any $b\in(-1,1)$ and any open interval $I\subseteq \R$,
the space
$$
L^{2;b}(I\to X):=L^2(I\to X;|y|^bdy)
$$
is endowed with the Hilbertian scalar product
$$
(U,V)_{L^{2;b}}=\intr|y|^b(U(y),V(y))_X~\!dy~\quad U,V\in L^{2;b}(I\to X)
$$
and corresponding norm $\|\cdot \|_{L^{2;b}}$ (notice that we avoid the longer notation $\|\cdot \|_{L^{2;b}(I\mapsto X)}$).

\medskip
Let $k\ge 0$ be an integer. We denote by ${\mathcal C}^k(I\to X)$  the space of curves $I\to X$ which are continuously differentiable up to the order $k$. If $U\in {\mathcal C}^k(I\to X)$, then 
$\partial^\ell_y U$ is the derivative of order $\ell=0,\dots,k$ (however, we will often write $\partial^2_{yy}$ instead of
$\partial^2_y$).
Further, ${\mathcal C}^\infty(I\to X)=\bigcap\limits_{k\ge 0} {\mathcal C}^k(I\to X)$. 

Accordingly with a commonly used notation, curves in ${\mathcal C}^{k,\sigma}(I\to X)\subset {\mathcal C}^k(I\to X)$ have H\"older continuous derivatives of order $k$. For our purposes, it is convenient to put
\begin{equation}
\label{eq:Holder_space}
\widetilde{\mathcal C}^{\alpha}(I\to X)=
\begin{cases}
{\mathcal C}^{\lfloor{\alpha}\rfloor,\alpha-\lfloor{\alpha}\rfloor}(I\to X)&\text{if $\alpha>0$ is not an integer}\\
{\mathcal C}^{\lfloor{\alpha}\rfloor-1,1}(I\to X)&\text{if $\alpha\ge 1$ is an integer.}
\end{cases}
\end{equation}
Also, for $U\in   \widetilde{\mathcal C}^{\alpha}(I\to X)$  we put
$$
\llbracket  U\rrbracket_{\widetilde{\mathcal C}^{\alpha}}=
\begin{cases}
\displaystyle{\sup_{y_1,y_2\in\R\atop y_1\neq y_2}\frac{\|\partial^{\isalpha}_yU(y_1)-\partial^{\isalpha}_yU(y_2)\|_X}{|y_1-y_2|^{\alpha-\isalpha}}}&\text{if $\alpha\notin \mathbb N$,}\\
\displaystyle{\sup_{y_1,y_2\in\R\atop y_1\neq y_2}\frac{\|\partial^{\alpha-1}_yU(y_1)-\partial^{\alpha-1}_yU(y_2)\|_X}{|y_1-y_2|}}
&\text{if $\alpha\in \mathbb N$.}
\end{cases}
$$

Notice that $\widetilde{\mathcal C}^{\alpha}(I\to X)
\subset {\mathcal C}^{\lfloor{\alpha}\rfloor}(I\to X)$ if and only if $\alpha$ is not an integer.}

\medskip
Let $k\in\mathbb N\cup\{\infty\}$. The spaces of even curves in $L^{2;b}(\R\to X), {\mathcal C}^k(\R\to X)$ are denoted by $L^{2;b}_\me(\R\to X)$, $ {\mathcal C}^k_\me(\R\to X)$, respectively,
and 
${\mathcal C}^k_{c,\me}(\R\to X)$ is the space of compactly supported functions in ${\mathcal C}^k_{\me}(\R\to X)$.

\medskip 
We simply write $L^{2;b}_\me(\R), {\mathcal C}^k_\me(\R), {\mathcal C}^\infty_{c,\me}(\R)$ instead of $L^{2;b}(\R\to \R), {\mathcal C}^k(\R\to \R), {\mathcal C}^\infty_{c,\me}(\R\to\R)$.

\normalsize 

\section{Spaces of real valued functions }
\label{S:spaces}

In this section, for any parameter $b\in (-1,1)$ and any integer $k\ge 0$ we introduce the 
Sobolev-type space $H^{k;b}_\me(\R)$, which is related to
the differential operators $(\Db_b+\lambda)^k$, $\lambda>0$.

The choice of working with even functions has been inspired by \cite{CM}. This strategy
is needed in case $b\neq 0$ to overcome some  technical difficulties produced by
 the singularity of the operator $\Db_b$ in (\ref{eq:Db}) at $y=0$.

In fact, as noticed in \cite{CM}, if $\psi\in {\mathcal C}^2_{\me}(\R)$, then $y^{-1}\partial_y\psi(y)=\partial^2_{yy}\psi(0)+o(1)$ as $y\to 0$, which implies 
$\Db_b\psi\in {\mathcal C}^0_{\me}(\R)$ More generally,
\begin{equation}
\label{eq:even0}
(\Db_b+\lambda)^m\psi \in {\mathcal C}^{k-2m}_{\me}(\R)\quad\text{for any integer $m\le k/2$ and any $\psi\in {\mathcal C}^{k}_\me(\R)$.}
\end{equation}
Our definition of $H^{k;b}_\me(\R)$ is based on induction procedure, starting from the lower order cases $k=1,2$.

\paragraph{First order}
For $\lambda>0$, we endow the weighted Hilbert space 
$$
H^{1;b}(\R):=H^{1}(\R;|y|^bdy)= \{\psi\in L^{2;b}(\R)~|~ \partial_y\psi \in L^{2;b}(\R)\}
$$
with the scalar product
$$
(\psi,\eta)_{\lambda,H^{1;b}}=\intr|y|^{b}(\partial_y\psi\partial_y\eta+\lambda\psi\eta)~\!dy
$$
and the corresponding norm $\|\psi\|_{\lambda,H^{1;b}}$. If $\lambda=1$ we drop it and simply write $(\psi,\eta)_{H^{1;b}}$ and $\|\psi\|_{H^{1;b}}$.
Clearly, the norms $\|\cdot \|_{\lambda,H^{1;b}}$ are equivalent for all $\lambda>0$
and moreover 
\begin{equation}
\label{eq:equivalent_k=1}
\|\psi(\cdot\sqrt{\lambda})\|_{\lambda,H^{1;b}_\me}^2=\lambda^{1-\frac{1+b}{2}}\|\psi(\cdot)\|_{H^{1;b}_\me}^2.
\end{equation}
 
\begin{Lemma}
\label{L:k=1}
\begin{itemize}
\item[$i)$] ${\mathcal C}^\infty_c(\R)$ is dense in $H^{1;b}(\R)$;
\item[$ii)$] $H^{1;b}(\R)\subset H^1_{\rm loc}(\R)$ if $b\in(-1,0]$ and $H^{1;b}(\R)\subset W^{1,p}_{\rm loc}(\R)$ for any $p\in[1,\frac{2}{1+{b}})$ if $b\in (0,1)$;
\item[$iii)$] $H^{1;b}(\R)\subset  {\mathcal C}^{0,\frac{1-b_+}{2}}_{\rm loc}(\R)$;
\item[$iv)$] There exists  $m_b>0$ such that $\|\psi\|^2_{H^{1;b}}\ge m_b|\psi(0)|^2$ for any $\psi\in H^{1;b}(\R)$.
\end{itemize}
\end{Lemma}

\proof 
For $i)$ see \cite{K}. 
The first part of $ii)$ is trivial; to prove the second one use H\"older's inequality.

If $b\le 0$ then  $ii)$ implies $iii)$ immediately. Assume $b\in(0,1)$ and take $\psi\in {\mathcal C}^\infty_{c}(\R)$. Since 
$$
\psi(y_2)-\psi({y_1})=\int\limits_{y_1}^{{y_2}} |t|^{\frac{-{b}}{2}}(|t|^{\frac{{b}}{2}}\partial_t\psi(t))\,dt
$$
for any $y_1, {y_2}\in\R$, then H\"older's inequality and the density result in $i)$ imply that  
\begin{equation}
\label{eq:Hol}
\begin{aligned}
|\psi(y_2)-\psi({y_1})|^2&\le
\frac{1}{1-{b}}\|\partial_y\psi\|^2_{L^{2;b}}~\!\big| y_2|y_2|^{-{b}}-{y_1}|{y_1}|^{-{b}}\big|\\&
\le
\frac{2^{1-b}}{1-{b}}\|\partial_y\psi\|^2_{L^{2;b}}~\!\big| y_2-y_1|^{1-{b}}
\end{aligned}
\end{equation}
for any $\psi \in H^{1;b}(\R)$,  $y_1,y_2\in\R$. Since ${\mathcal C}^\infty_c(\R)$ is dense in $H^{1;b}(\R)$ and since $\psi$ was arbitrarily chosen in $H^{1;b}(\R)$, the inclusion in $iii)$ 
easily follows.

Lastly, given $\psi\in H^{1;b}(\R)$ we use (\ref{eq:Hol}) to get the existence of a constant  $c>0$ depending only on $b$ such that
$$
c|y|^b|\psi(0)|^2\le |y|\|\partial_y\psi\|^2_{L^{2;b}}+|y|^b|\psi(y)|^2
$$
for any $y\in\R$. Then $iii)$ follows via integration over $(0,1)$.
\QED

\begin{Remark}
\label{R:c}
It follows from Theorem \ref{T:Aexistence} in Section \ref{S:Bessel} that the best constant in $iv)$
is 
$$
m_{b}=2^{1+b}\Gamma\Big(\frac{1+b}{2}\Big)\Gamma\Big(\frac{1-b}{2}\Big)^{-1},
$$
and it is achieved by the function $\psi_s$, see (\ref{eq:psi_notation}), for $s=\frac{1-b}{2}$.
\end{Remark}

We will be mainly concerned with $H^{1;b}_\me(\R)$, the subspace of even functions in $H^{1;b}(\R)$.  For future convenience, we notice that the proof of Lemma \ref{L:k=1} gives
\begin{equation}
\label{eq:loc_Holder}
|\psi(y_2)-\psi({y_1})|^2\le \frac{1}{1-{b}}\|\partial_y\psi\|^2_{L^{2;b}}~\!\big| |y_2|^{1-{b}}-|{y_1}|^{1-{b}}\big|
\end{equation}
for any $\psi\in H^{1;b}_\me(\R)$, $y_1, y_2\in\R$.

\paragraph{Second order}
If $\psi\in H^{1;b}_\me(\R)$ then $|y|^b\partial_y\psi\in L^{2;-b}(\R)\subset L^1_{\rm loc}(\R)$. We put
$$
H^{2;b}_\me(\R)=\big\{\psi\in H^{1;b}_\me(\R)~|~ |y|^b\partial_y\psi\in H^{1;-b}(\R) \big\}.
$$
Let $\psi\in {\mathcal C}^2_{c,\me}(\R)$. Then $\partial_y(|y|^b\partial_y\psi)= -|y|^{b}\Db_b\psi$, which implies $\psi \in H^{2;b}_\me(\R)$ by (\ref{eq:even0}).
We 
extend the pointwise defined operator
$\Db_b$ to $H^{2;b}_\me(\R)$ by putting
$$
\Db_b\psi:=-|y|^{-b}\partial_y(|y|^b\partial_y\psi)\qquad\text{for $\psi\in H^{2;b}_\me(\R)$,}
$$
so that $\Db_b: H^{2;b}_\me(\R)\to L^{2;b}_\me(\R)$.

\begin{Lemma}
\label{L:parts2_psi}
Let $\psi\in H^{2;b}_\me(\R)$. Then
\begin{eqnarray}
\label{eq:ii}
&(\Db_b\psi, \eta)_{L^{2;b}}~\!=(\partial_y\psi ,\partial_y\eta)_{L^{2;b}}\quad &\text{for any $\eta\in H^{1;b}_\me(\R)$;}\\
\label{eq:iii} 
&(\Db_b\psi,\eta)_{L^{2;b}}=(~\!\psi~\!,\Db_b\eta)_{L^{2;b}}\quad &\text{for any $\eta\in H^{2;b}_\me(\R)$.}
\end{eqnarray}
\end{Lemma}

\proof 
Let $\eta\in {\mathcal C}^\infty_{c,\me}(\R)$. We can use integration by parts to compute
$$
\intr |y|^b(\Db_b\psi)\eta~\!dy=- \intr \partial_y(|y|^b\partial_y\psi)\eta~\!dy=
\intr |y|^b\partial_y\psi~\!\partial_y\eta~\!dy.
$$
Thus $i)$ follows, thanks to the density result in Lemma \ref{L:k=1}. Clearly $ii)$ is an immediate consequence of $i)$.\QED

It remains to introduce a Hilbertian structure on $H^{2;b}_\me(\R)$. Given $\lambda>0$, we put
$$
(\psi,\eta)_{\lambda,H^{2;b}_\me}=((\Db_b+\lambda)\psi,(\Db_b+\lambda)\eta)_{L^{2;b}}~,\quad 
\|\psi\|_{\lambda,H^{2;b}_\me}= \|(\Db_b+\lambda)\psi\|_{L^{2;b}}.
$$
If $\lambda=1$ we drop it and simply write $(\psi,\eta)_{H^{2;b}_\me}$ and $\|\psi\|_{H^{2;b}_\me}$. 
Notice that
\begin{equation}
\label{eq:rescaling0}
(\Db_b+\lambda)\psi(\cdot~\!\sqrt{\lambda})=\lambda~\!\big[(\Db_b+1)\psi\big](\cdot~\!\sqrt\lambda)~,
\end{equation}
which implies
\begin{equation}
\label{eq:rescaling2}
\|\psi(\cdot\sqrt{\lambda})\|_{\lambda,H^{2;b}_\me}^2=\lambda^{2-\frac{1+b}{2}}\|\psi(\cdot)\|_{H^{2;b}_\me}^2\quad
\text{for any $\psi\in H^{2;b}_\me(\R)$.}
\end{equation}

\medskip

\begin{Lemma}
\label{L:induction_2}
Let $\lambda>0$, $\psi \in H^{2;b}_\me(\R)$. Then  
$$\|\psi\|^2_{\lambda,H^{2;b}_\me}\ge  \lambda\|\psi\|^2_{\lambda,H^{1;b}}.
$$
Therefore, $H^{2;b}_\me(\R)$ is  a Hilbert space, and is continuously embedded in $H^{1;b}_\me(\R)$.
\end{Lemma}

\proof
Thanks to (\ref{eq:rescaling2}) we can assume that $\lambda=1$.
By Lemma \ref{L:parts2_psi} with $\eta=\psi$ we have
$(\Db_b\psi,\psi)_{L^{2;b}}= \|\partial_y\psi\|^2_{L^{2;b}}$. Thus
$$
\begin{aligned}
\intr|y|^b|(\Db_b+1)\psi|^2~\!dy&= \intr|y|^b|\Db_b\psi|^2~\!dy + 2\intr|y|^b(\Db_b\psi)\psi~\!dy+\!\!\!\intr|y|^b|\psi|^2~\!dy\\
&\ge 
2\intr |y|^b|\partial_y\psi|^2~\!dy+\!\!\!\intr|y|^b|\psi|^2~\!dy
\end{aligned}
$$
which implies $\|\psi\|^2_{H^{2;b}_\me}\ge \|\psi\|^2_{H^{1;b}}$. The conclusion of the proof is standard.
\QED

\paragraph{Higher order}
If $k>2$ and $\lambda>0$ we use induction to define
\begin{gather*}
H^{k;b}_\me(\R)=\Big\{\psi\in H^{k-1;b}_\me(\R)~|~ \Db_b\psi\in H^{k-2,b}_\me(\R)~\Big\}\\
(\psi,\eta)_{\lambda,H^{k;b}_\me}=((\Db_b+\lambda)\psi,(\Db_b+\lambda)\eta)_{\lambda,H^{k-2;b}_\me}~,\quad 
\|\psi\|_{\lambda,H^{k;b}_\me}= \|(\Db_b+\lambda)\psi\|_{\lambda,H^{k-2;b}_\me}.
\nonumber
\end{gather*}
As usual, if $\lambda=1$ we drop it and simply write $(\psi,\eta)_{H^{k;b}_\me}$ and $\|\psi\|_{H^{k;b}_\me}$.

Notice that ${\mathcal C}^k_{c,\me}(\R)\subset H^{k;b}_\me(\R)$ by (\ref{eq:even0}). In the next lemma we collect the main 
properties of the spaces $H^{k;b}_\me(\R)$ for $k\ge 1$. In particular it implies
that $\|\cdot \|_{\lambda,H^{k;b}_\me}$, for different $\lambda$'s, define the same Hilbertian structure on  $H^{k;b}_\me(\R)$.
We omit its easy proof, which is based on  previous results and induction.

\begin{Lemma}
\label{L:all}
Let $k\ge 1$, $b\in(-1,1)$, $\psi \in H^{k;b}_\me(\R)$ and $\lambda>0$. The following facts hold.
\begin{itemize}
\item[$i)$] ~$\displaystyle{
\|\psi\|_{\lambda,H^{k;b}_\me}^2 =
\begin{cases}
\|\partial_y(\Db_b+\lambda)^{\frac{k-1}{2}}\psi)\|^2_{L^{2;b}}+\lambda\|(\Db_b+\lambda)^{\frac{k-1}{2}}\psi\|^2_{L^{2;b}}&\text{if $k$ is  odd}\\
\|(\Db_b+\lambda)^{\frac{k}{2}}\psi\|^2_{L^{2;b}}&\text{if $k$ is even~\!;}
\end{cases}
}$
\item[$ii)$] 
~$\displaystyle{
\|\psi(\cdot\sqrt{\lambda})\|_{\lambda,H^{k;b}_\me}^2=\lambda^{k-\frac{1+b}{2}}\|\psi(\cdot)\|_{H^{k;b}_\me}^2~\!;
}$
\item[$iii)$]~$(\Db_b+\lambda)^m\psi\in H^{k-2m;b}_\me(\R)$ for any positive integer $m<k/2$;
\item[$iv)$] 
~$\displaystyle{
\|\psi\|^2_{\lambda,H^{k;b}_\me}\ge \lambda^{k-j}\|\psi\|^2_{\lambda,H^{j;b}}}\ge  \lambda^{k}\|\psi\|^2_{L^{2;b}}$ for any $j=1,\dots,k$;
\item[$v)$] 
~$\displaystyle{
\|\psi\|_{\lambda,H^{k;b}_\me}^2\ge m_b \lambda^{k-\frac{1+b}{2}}|\psi(0)|^2
}$, where $m_b$ is the constant in Lemma \ref{L:k=1}.
\end{itemize}
\end{Lemma}

We now establish some integration by parts formulae. It suffices to take $\lambda=1$.

\begin{Lemma}
\label{L:parts_psi}
Let $k\ge 2$, $\psi\in {H^{2(k-1);b}_\me(\R)} $, $\eta\in H^{k;b}_\me(\R)$.  Then 
$$(\psi,\eta)_{H^{k;b}_\me}=((\Db_b+{1})^{k-1}\psi,(\Db_b+{1})\eta)_{L^{2;b}}.$$
\end{Lemma}

\proof 
Notice that ${H^{2(k-1);b}_\me(\R)} \subset H^{k;b}_\me(\R)$. 

If $k=2$, the equality in the lemma holds by definition.

If $k=2m\ge 4$ is even, we use (\ref{eq:iii}) with  $(\Db_b+{1})^m\psi\in H^{2(m-1);b}_{\me}(\R)$ instead of $\psi$ and 
$(\Db_b+{1})^{m-1}\eta\in H^{2;b}_\me(\R)$ instead of $\eta$ to get
$$
(\psi,\eta)_{H^{2m;b}_\me}=((\Db_b+{1})^{m}\psi,(\Db_b+{1})^{m}\eta)_{L^{2;b}}=
((\Db_b+{1})^{m+1}\psi,(\Db_b+{1})^{m-1}\eta)_{L^{2;b}}.
$$
If $m=2$ we are done. Otherwise, repeat the same procedure $m-1$ times to get 
\begin{equation}
\label{eq:even}
(\psi,\eta)_{H^{2m;b}_\me}= 
((\Db_b+{1})^{2m-1}\psi,(\Db_b+{1})\eta)_{L^{2;b}}~,
\end{equation}
which concludes the proof in the even case.

If $k=2m+1\ge 3$ is odd 
we apply 
(\ref{eq:ii}) with $(\Db_b+{1})^{m}\psi\in H^{2m;b}_{\me}(\R)$ instead of $\psi$ and $(\Db_b+{1})^{m}\eta\in H^{1;b}_\me(\R)$ instead of $\eta$ to get
$$
(~\!\partial_y((\Db_b+{1})^{m}\psi),\partial_y((\Db_b+{1})^{m}\eta))_{L^{2;b}}
=
(\Db_b(\Db_b+{1})^{m}\psi, (\Db_b+{1})^{m}\eta)_{L^{2;b}}.
$$
It follows that
$$
\begin{aligned}
(\psi,\eta)_{H^{k;b}_\me}&= (~\!\partial_y((\Db_b+{1})^{m}\psi),\partial_y((\Db_b+{1})^m \eta))_{L^{2;b}}
+((\Db_b+{1})^{m}\psi,(\Db_b+{1})^m \eta)_{L^{2;b}}\\
&=
((\Db_b+{1})^{m+1}\psi,(\Db_b+{1})^{m}\eta)_{L^{2;b}}=
((\Db_b+{1})\psi,\eta)_{H^{2m;b}_\me}.
\end{aligned}
$$
To conclude the proof, use (\ref{eq:even}) with $\psi$ replaced by $(\Db_b+{1})\psi$.
\QED

\begin{Remark} 
It is well known that smooth, compactly supported functions are dense in $H^k(\R)$ for any $k>0$. Recall that  
${\mathcal C}^\infty_{c}(\R)$ is dense in $ H^{1;b}(\R)$ for any $b\in(-1,1)$ by \cite{K}.
If would be of interest to prove the density of ${\mathcal C}^\infty_{c,\me}(\R)$ in $ H^{k;b}_\me(\R)$ in case $b\neq 0$, $k>1$.
\end{Remark}

\section{Bessel functions and related issues}
 \label{S:Bessel}
 
The basic properties of the Bessel function $K_\alpha$ can be found for instance \cite[Sections 8.4, 8.5]{Grad}. 
For any $\alpha\in \R$ the standard modified Bessel function of the second kind $K_\alpha=K_{-\alpha}$ solves
$$
\partial^2_{yy} K_\alpha(y) +y^{-1}\partial_y K_\alpha(y)-(1+\alpha^2y^{-2})K_\alpha(y)=0\quad
\text{on $\R_+$}
$$ 
and decays exponentially as $y\to+\infty$. If $\alpha\neq 0$ then
$$
K_\alpha(y)=2^{|\alpha|-1}\Gamma(|\alpha|)y^{-|\alpha|} +o(y^{-|\alpha|} )\quad
\text{as $y\to 0^+$.}
$$
Bessel functions of different orders are related by the formulae
$$
\partial_y(y^\alpha K_\alpha(y))=-y^\alpha K_{\alpha-1}(y)~,\qquad 
K_{\alpha}(y)-K_{{\alpha}-2}(y)=2({\alpha}-1)y^{-1}K_{{\alpha}-1}(y).
$$

Next, for $s>0$ and $\lambda>0$ we put 
\begin{equation}
\label{eq:psilambda}
\psi_{s,\lambda}(y):=\psi_s(\sqrt{\lambda}~\!y)=c_s (\sqrt{\lambda}|y|)^s K_s(\sqrt{\lambda}|y|)~\!,
\end{equation}
see  (\ref{eq:psi_notation}). 
Notice that 
$$
\psi_{s,\lambda}\in {\mathcal C}^0_\me(\R)~,\quad \psi_{s,\lambda}(0)=1~,\quad  \psi_{s,\lambda}\in {\mathcal C}^\infty(\R_+)~\!,
$$
and $\psi_{s,\lambda}$ decays exponentially at infinity together with its derivatives of any order. Further, (\ref{eq:rescaling0}) 
readily implies
\begin{equation}
\label{eq:rescaling100}
(\Db_b+\lambda)^m\psi_{s,\lambda}(y)=\lambda^m~\!\big[(\Db_b+1)^m\psi_{s}\big](\sqrt\lambda~\!y)
\end{equation}
for any $y\neq 0$ and any integer $m\ge 1$.

\begin{Lemma}
\label{L:old}
Let $s>0$ be non-integer and put ${\mathfrak b}=1-2(s-\lfloor{s}\rfloor)$. Then 
$\psi_s$ solves the following
differential equations on $\R_+$:
\begin{itemize}
\item[$i)$] 
$
\partial_y\psi_s(y)=
\begin{cases}
-d_s~\! y^{2s-1}\psi_{1-s}(y)&\text{if ~$0<s<1$,}\\
-\frac{1}{2(s-1)}~\!y\psi_{s-1}(y)&\text{if ~$s>1$;}
\end{cases}$
\item[$ii)$]
$
-\partial^2_{yy}\psi_s(y)+ \psi_s(y)=
\begin{cases}
d_s(2s-1)y^{2(s-1)}\psi_{1-s}(y)& \text{if ~$0<s<1$,}\\
\frac{2s-1}{2(s-1)} \psi_{s-1}& \text{if ~$s>1$;}
\end{cases}
$
\item[$iii)$]
$(\Db_{\mathfrak b}+1)^{{\lceil{s}\rceil}}\psi_{s}=0$;
\item[$iv)$] If $s>1$ then for any  $m=1,\dots, \lfloor{s}\rfloor$
\begin{equation}
\label{eq:Abm}
(\Db_{\mathfrak b}+1)^m\psi_{s}=  
\frac{d_s}{d_{s-m}}~\!\psi_{{s}-m}= \frac{\lfloor{s}\rfloor!}{\lfloor{s-m}\rfloor!}~\!\frac{\Gamma(s-m)}{\Gamma(s)}~\!\psi_{{s}-m}.
\end{equation}
\end{itemize}
\end{Lemma}

\proof
Let $s\in(0,1)$. By the properties of the Bessel functions we get 
$$
\partial_y\psi_{s}(y)= -c_{s} y^{s} K_{1-s}(y)=  -c_{s} y^{2s-1} (y^{1-s}K_{1-s}(y))=-d_sy^{2s-1}\psi_{1-s}(y).
$$ 
This gives the first equality in $i)$.
Now we notice that we can compute $\partial_y\psi_{1-s}$ via {the first equality in $i)$,
where $s$ is replaced by $1-s$. The proofs of $ii)$, $iii)$ readily follow.
This completes the proof in this case. }

Now let $s>1$. We compute
$$
\partial_y\psi_{s}(y)= c_{s} \partial_y(y^{s} K_{s}(y))= -c_{s}y (y^{s-1} K_{{s}-1}(y))=-\frac{c_s}{c_{s-1}}\psi_{s-1}(y)
$$
which gives the second equality in $i)$. 
Also, we get
$$
\begin{aligned}
\partial^2_{yy}\psi_{s}(y)&= -c_{s} \partial_y(y^{s} K_{{s}-1}(y))=-c_sy^{s} (- K_{{s}-2}(y)+y^{-1}K_{{s}-1}(y))\\ &
=c_sy^s\big((1-2s)y^{-1}K_{s-1}(y)+K_s(y)\big)
\end{aligned}
$$
by the recurrence formula for $K_s$. Hence 
$$\partial^2_{yy}\psi_{s}(y)= \frac{c_s}{c_{s-1}}(1-2s)\psi_{s-1}+\psi_s(y),
$$
which gives $ii)$ for $s>1$. To prove $iv)$ we notice  that this last equality  implies
$$
\begin{aligned}
(\Db_{\mathfrak b} +1)\psi_s&=
-\partial^2_{yy}\psi_{s}-(1-2s+2\lfloor{s}\rfloor)\partial_{y}\psi_{s}+\psi_s=
\frac{\lfloor{s}\rfloor}{s-1}\psi_{s-1}=
\frac{d_s}{d_{s-1}}\psi_{s-1}~\!.
\end{aligned}
$$
Thus (\ref{eq:Abm}) holds for $m=1$. To conclude the proof of $iv)$ repeat the same argument a finite number of times.

It remains to prove $iii)$ in this case. We use $iv)$ with $m= \lfloor{s}\rfloor$ and then $iii)$ with $s$ replaced by $s-\lfloor{s}\rfloor\in(0,1)$ to get
$$
(\Db_{\mathfrak b} +1)^{\lceil{s}\rceil}\psi_s 
=\frac{d_s}{d_{1-\lfloor{s}\rfloor}}  (\Db_{\mathfrak b} 
+1)\psi_{s-\lfloor{s}\rfloor}=0~\!.
$$
The lemma is completely proved.
\QED

\begin{Remark}
\label{R:max_psi}
Since $K_s>0$ on $\R_+$, from $i)$ in Lemma \ref{L:old} it readily follows that the positive function $\psi_s$ achieves its maximum at
the origin.
\end{Remark}

The next theorem contains our main result on the functions $\psi_s$ (recall our non-standard definition of  H\"older spaces in (\ref{eq:Holder_space})).

\begin{Theorem}
\label{T:Aexistence} 
Let $s>0$ be non-integer, put ${\mathfrak b}=1-2(s-\lfloor{s}\rfloor)$ and let 
$\lambda>0$. Then
\begin{gather}
\label{eq:regularity}
\psi_{s,\lambda}\in {H}^{\lceil{s}\rceil;{\mathfrak b}}_\me(\R)\cap 
{\widetilde{\mathcal C}^{2s}(\R)}~;\\
\label{eq:operators}
\lim_{y\to 0^+}y^{\mathfrak b}\partial_y\big((\Db_{\mathfrak 
b}+\lambda)^{\lfloor{s}\rfloor} \psi_{s,\lambda})= 
- ~\!d_s\lambda^{s}
\end{gather}
where $d_s$ is the constant  in (\ref{eq:isometry}). 
Moreover $\psi_{s,\lambda}$ satisfies
\begin{equation}
\label{eq:equation_psi}
(\psi_{s,\lambda},\eta)_{\lambda,H^{\lceil{s}\rceil;{\mathfrak 
b}}_\me}=~\!2d_s~\!\lambda^s~\!\eta(0)\quad\text{for any} \quad\eta\in 
{H}^{\lceil{s}\rceil;{\mathfrak b}}_\me(\R).
\end{equation}
Finally, $\psi_{s,\lambda}$  admits the following variational characterization, 
\begin{gather}
\label{eq:minimum}
\|\psi_{s,\lambda}\|^2_{\lambda,H^{\lceil{s}\rceil;{\mathfrak b}}_\me}=  
\inf_{\eta\in {H}^{\lceil{s}\rceil;{\mathfrak b}}_\me(\R)\atop \eta(0)=1} 
\|~\!\eta~\!\|^2_{\lambda,H^{\lceil{s}\rceil;{\mathfrak b}}_\me} =2d_s 
\lambda^s~\!.
\end{gather}
\end{Theorem}

\proof
Thanks to (\ref{eq:rescaling100}), we  assume that $\lambda=1$.
We divide the proof in two steps.

\medskip
\noindent{\bf Step 1.} Let $\lfloor{s}\rfloor=0$. Then ${\mathfrak b}=1-2s$ and 
\begin{equation}
\label{eq:psi-alpha1}
\partial_y\psi_{s}(y)=-d_s~\! y^{-{\mathfrak b}}\psi_{1-s}(y)=- 
d_s~\!y^{-{\mathfrak b}}+o(y^{-{\mathfrak b}})\quad\text{as $y\to 0^+$,}
\end{equation}
which proves (\ref{eq:operators}).
Since in addition $\psi_s$ decays exponentially at infinity, from (\ref{eq:psi-alpha1}) we first infer that $\psi_{s}\in H^{1;1-2{s}}_\me(\R)$. 

\medskip

To prove that $\psi_s\in \widetilde{\mathcal C}^{2s}(\R)$ we fix two points $y_1,y_2\in\R$. By the symmetry of $\psi_s$, we can assume that $y_1,y_2\ge 0$. 

Let $0<2s\le 1$.
For $y>0$ we have $|\partial_y\psi_s(y)|=d_s y^{2s-1}\psi_{1-s}(y)\le d_s y^{2s-1}$. Thus
$\psi_s\in \widetilde{\mathcal C}^{2s}(\R)$ follows from
\begin{equation*}
|\psi_s(y_1)-\psi_s(y_2)|\le d_s\Big|\int_{y_1}^{y_2}y^{2s-1}dy\Big|= \frac{d_s}{2s}|y_1^{2s}-y_2^{2s}|
\le \frac{d_s}{2s}|y_1-y_2|^{2s}.
\end{equation*}

If  $1<2s<2$ we use $ii)$ in Lemma \ref{L:old} to estimate
$$
|\partial^2_{yy}\psi_s(y)|=|\psi_s(y)-d_s(2s-1)y^{2(s-1)}\psi_{1-s}(y)|
\le  1+cy^{2(s-1)}
$$
for $y>0$. Using integration as before, we plainly get 
$$
|\partial_y\psi_s(y_1)-\partial_y\psi_s(y_2)|\le 
|y_1-y_2|+c|y_1-y_2|^{2s-1}.
$$
Since $\partial_y\psi_s$ decays exponentially at infinity, we infer that 
there exists a constant $c>0$ depending only on $s$, such that
\begin{equation*}
|\partial_y\psi_s(y_1)-\partial_y\psi_s(y_2)|\le c|y_1-y_2|^{2s-1}
\end{equation*}
which, in turns concludes the proof of (\ref{eq:regularity}).

Next, by $iii)$ in Lemma \ref{L:old} we have that 
\begin{equation}
\label{eq:psi_eq}
\partial_y(y^{\mathfrak b}\partial_y\psi_s)= y^{\mathfrak b}\psi_s\quad 
\text{on $\R_+$.}
\end{equation}
We test (\ref{eq:psi_eq}) with an arbitrary $\eta\in {\mathcal C}^\infty_{c,\me}(\R)$. Taking (\ref{eq:psi-alpha1}) into account we obtain
$$
\int\limits_{0}^\infty y^{\mathfrak b}\psi_s\eta~\!dy= 
\int\limits_0^\infty\partial_y(y^{\mathfrak 
b}\partial_y\psi_s)\eta~\!dy=d_s\eta(0) -\int\limits_0^\infty 
y^{\mathfrak b}\partial_y\psi_s\partial_y\eta~\!dy.
$$
By the evenness of $\psi_s$ and $\eta$, this implies that
$(\psi_s,\eta)_{H^{1;{\mathfrak b}}_\me}=2d_s\eta(0)$. Thus 
(\ref{eq:equation_psi}) holds in case $\lfloor{s}\rfloor=0$, thanks to
the density result in Lemma \ref{L:k=1}.

From (\ref{eq:equation_psi}) it follows that 
$(\psi_s,\eta-\psi_s)_{H^{1;{\mathfrak b}}_\me}=0$ for any $\eta\in 
H^{1;{\mathfrak b}}_\me(\R)$ such that $\eta(0)=1$.
Thus, $\psi_s$ is the minimal distance projection of $0$ on the hyperplane 
$\{\eta(0)=1\}\subset H^{1;{\mathfrak b}}_\me(\R)$, that is, $\psi_s$ is the 
unique solution to the minimization problem in (\ref{eq:minimum}).
This completes the proof in the case $s\in(0,1)$.

\medskip
\noindent{\bf Step 2:} Let $\lfloor{s}\rfloor\ge 1$. 
Thanks to $i)$ in Lemma \ref{L:old} we see that 
$$
\partial_y\psi_s(y)=-\frac{1}{2(s-1)}y~\psi_{s-1}(y)= -\frac{1}{2(s-1)}y+o(y)\quad\text{as $y\to 0^+$,}
$$
hence $\psi_s\in {\mathcal C}^{2}(\R)$. Next, as in case $2s\in(1,2)$ we use 
$ii)$ in Lemma \ref{L:old} to infer that $\partial^2_{yy}\psi_s$ has the same regularity as $\psi_{s-1}$. 
If $s\in(1,2)$ we  
obtain $\psi_s\in \widetilde{\mathcal C}^{2s}(\R)$ by Step 1; if $s>2$ one can use a bootstrap argument to prove that 
$\psi_{s}\in \widetilde{\mathcal C}^{2s}(\R)$. By the decaying of $\psi_{s}$ at infinity we also infer that
$$
\psi_{s}\in H^{2\lfloor{s}\rfloor;{\mathfrak b}}_\me(\R)\subset 
H^{\lceil{s}\rceil;{\mathfrak b}}_\me(\R),$$
which concludes the proof of (\ref{eq:regularity}).

\medskip
To prove (\ref{eq:operators}) it suffices to notice that (\ref{eq:Abm}) and Step 1 give
\begin{equation*}
\lim_{y\to 0^+} y^{\mathfrak b}\partial_y\big((\Db_{\mathfrak 
b}+1)^{\lfloor{s}\rfloor} 
\psi_{s})=\frac{d_s}{d_{s-\lfloor{s}\rfloor}}\lim_{y\to 0^+} y^{\mathfrak b} 
\partial_y\psi_{s-\lfloor{s}\rfloor}=
-d_s.
\end{equation*}

We now prove  (\ref{eq:equation_psi}).   Take any $\eta\in 
H^{\lceil{s}\rceil;{\mathfrak b}}_\me(\R)$. We apply Lemma \ref{L:parts_psi} 
with $k=\lceil{s}\rceil$ and $\psi=\psi_s$ to obtain
$$
(\psi_s,\eta)_{H^{\lceil{s}\rceil;{\mathfrak 
b}}_\me}=((\Db_{\mathfrak 
b}+1)^{\lfloor{s}\rfloor}\psi_s,(\Db_{\mathfrak b}+1)\eta)_{L^{2;{\mathfrak 
b}}}.
$$
Therefore, (\ref{eq:Abm}),  (\ref{eq:ii}) and Step 1 with $s$ replaced by $s-\lfloor{s}\rfloor\in (0,1)$ give
$$
(\psi_s,\eta)_{H^{\lceil{s}\rceil;{\mathfrak 
b}}_\me}=\frac{d_s}{d_{s-\lfloor{s}\rfloor}}(\psi_{s-\lfloor{s}\rfloor},
(\Db_{\mathfrak b}+1)\eta)_{L^{2;{\mathfrak b}}}=
\frac{d_s}{d_{s-\lfloor{s}\rfloor}}(\psi_{s-\lfloor{s}\rfloor},\eta)_{H^{1;{
\mathfrak b}}_\me}=2d_s\eta(0),
$$
and (\ref{eq:equation_psi}) follows.  For (\ref{eq:minimum}) argue as in Step 1.
\QED

\begin{Remark}
The recurrence formulae (\ref{eq:Abm}) plainly imply the  identities
$$
\begin{aligned}
(\Db_{\mathfrak b}+1)^m\psi_s(0)&=\frac{d_s}{d_{s-m}}~,\qquad 
m=1,\dots,\lfloor{s}\rfloor
\\
 \lim_{y\to 0^+}y^{-1}\partial_y\big((\Db_{\mathfrak b}+1)^{m} \psi_s)&= 
-\frac{d_s}{d_{s-m}}~\!\frac{1}{2(s-m-1)}~,\qquad m=0,\dots,\lfloor{s}\rfloor-1
\end{aligned}
$$
\end{Remark}

We conclude this section with a corollary of Theorem \ref{T:Aexistence}, which
might be of independent interest. 

\begin{Corollary}
Let {$\lfloor{s}\rfloor$ be even}. Then the following virial-type formulae hold:
\begin{gather*}
\intr |y|^{\mathfrak 
b}|(\Db_{\mathfrak 
b}+1)^\frac{\lfloor{s}\rfloor}{2}\psi_s|^2=\frac{s}{\lceil{s}\rceil} 
~\!2d_s\!~,\\
\intr |y|^{\mathfrak 
b}|\partial_y((\Db_{\mathfrak 
b}+1)^\frac{\lfloor{s}\rfloor}{2}\psi_s)|^2=\frac{\lceil{s}\rceil-s}{\lceil{s}\rceil}~\!2d_s~\!. 
\end{gather*}
\end{Corollary}

\proof
We use (\ref{eq:minimum}) with $s$ replaced by $s+1$ and then (\ref{eq:Abm}) to get
$$
\begin{aligned}
2d_{s+1}&=\|\psi_{s+1}\|_{H^{2+\lfloor{s}\rfloor;{\mathfrak 
b}}_\me}^2=\intr|y|^{\mathfrak 
b}|(\Db_{\mathfrak b}+1)^{\frac{\lfloor{s}\rfloor}{2}+1}\psi_{s+1} |^2~\!dy\\&
=
\frac{d^2_{s+1}}{d^2_s} \intr|y|^{\mathfrak 
b}|(\Db_{\mathfrak b}+1)^{\frac{\lfloor{s}\rfloor}{2}}\psi_s|^2~\!dy,
\end{aligned}
$$
and the first equality follows. For the second one, recall that 
$2d_s=\|\psi_{s}\|_{H^{\lceil{s}\rceil;{\mathfrak b}}_\me}^2$.
\QED

\section{Spaces of curves in $\mathcal H$; proof of the main results}
\label{SS:U_space}

We start this section by studying the (unbounded) operators
$\Lb_b^{\!k}U=(\Db_b+\mathcal L)^kU$ on  $L^{2;b}_\me(\R\to \mathcal H)$, for any 
$b\in(-1,1)$ and any integer $k\ge 0$. 

\medskip

Any function $U\in L^{2;b}(\R\to\mathcal H)$ can be decomposed as follows,
$$
U(y)=\sum_{j=1}^\infty U_j(y)\f_j,
$$
where $U_j=(U,\f_j)_{\!\mathcal H}\in L^{2;b}_\me(\R)$ for any  $j\ge 1$, and
$$
\begin{aligned}
\|U\|^2_{L^{2;b}}
= \sum_{j=1}^\infty\intr|y|^b|U_j|^2~\!dy
= \sum_{j=1}^\infty\|U_j\|_{L^{2;b}}^2.
\end{aligned}
$$
Recall that 
$$\Lb_bU=(\Db_b+\mathcal L)U= -\partial^2_{yy}U-by^{-1}\partial_yU+\mathcal LU~,\quad \mathcal L\f_j=\lambda_j\f_j
$$ and that 
we are assuming $\lambda_j\ge \lambda_1>0$.
Thus, at least formally we have
$$
\Lb_b^{\!k}U= \sum_{j=1}^\infty[(\Db_b+\lambda_j)^k U_j]~\!\f_j.
$$
We define
$$
H^{k;b}_{\mathcal L,\me}(\R\to\mathcal H)=\Big\{U\in L^{2;b}_\me(\R\to\mathcal H)~\big|~
U_j=(U,\f_j)_{\!\mathcal H}\in H^{k;b}_\me(\R) ~\text{and}~\|U\|_{H^{k;b}_{{\mathcal L},\me}}<\infty
\!~\Big\},
$$
where
$$
\|U\|^2_{H^{k;b}_{\mathcal L,\me}} 
:=\sum_{j=1}^\infty~\|U_j\|_{\lambda_j,H^{k;b}_\me}^2
$$
(we recall that $\|\cdot\|_{\lambda_j,H^{k;b}_\me}$ are equivalent norms in the space $H^{k;b}_\me(\R)$, see Section \ref{S:spaces}).
Thanks to Lemma \ref{L:all}, it is easily checked that $H^{k;b}_{\mathcal L,\me}(\R\to\mathcal H)$ is a Hilbert space with scalar product
$$
(U,V)_{\!H^{k;b}_{\mathcal L,\me}}=\sum_{j=1}^\infty~(U_j,V_j)_{\lambda_j,H^{k;b}_\me}.
$$

For future convenience we provide another definition of the spaces $H^{k;b}_{\mathcal L,\me}(\R\to\mathcal H)$. 
Consider the standard weighted Sobolev space
$$
H^{1;b}(\R\to\mathcal H) := H^{1}(\R\to\mathcal H;|y|^bdy)=\{U\in L^{2;b}(\R\to\mathcal H)~|~ \partial_yU \in L^{2;b}(\R\to\mathcal H)\},
$$
and denote by $H^{1;b}_\me(\R\to\mathcal H)$ the space of even curves in $H^{1;b}(\R\to\mathcal H)$. Then we let
$$
H^{2;b}_\me(\R\to\mathcal H)=\{U\in H^{1;b}_\me(\R\to\mathcal H)~|~|y|^b\partial_y U\in H^{1;-b}(\R\to\mathcal H)\} 
$$
so that 
$$
\Db_bU:=-|y|^{-b}\partial_y(|y|^b\partial_yU)\in L^{2;b}_\me(\R\to\mathcal H)\quad \text{for any $U\in H^{2;b}_\me(\R\to\mathcal H)$.}
$$
Finally, for $k\ge 3$ we use induction to define
$$
H^{k;b}_\me(\R\to\mathcal H)=\{U\in H^{k-1;b}_\me(\R\to\mathcal H)~|~\Db_b U\in H^{k-2;b}_\me(\R\to\mathcal H)\}.
$$

The proof of the next lemma is simple but boring, and we omit it.

\begin{Lemma}
\label{L:same_spaces}
Let $k\ge 1$ be an integer, $b\in(-1,1)$. Then
\begin{equation*}
H^{k;b}_{\mathcal L,\me}(\R\to\mathcal H)= H^{k;b}_{\me}(\R\to\mathcal H)\cap L^{2;b}(\R\to \mathcal H^k_\mathcal L).
\end{equation*}
\end{Lemma}

The next lemma will be useful for the proof of our main results.

\begin{Lemma}
\label{def:Lspace}
\begin{itemize}
\item[$i)$] If $U \in H^{k;b}_{\mathcal L,\me}(\R\to\mathcal H)$ then the following facts hold,
\begin{eqnarray}
\nonumber
\|U\|_{H^{k;b}_{\mathcal L,\me}}^2&=&
\begin{cases} 
\displaystyle{\| \partial_y(\Lb_b^{\frac{k-1}{2}}U)\|^2_{L^{2;b}}+\| {\mathcal L}^{\frac12}(\Lb_b^{\frac{k-1}{2}}U)\|^2_{L^{2;b}}}
&\text{if $k$ is odd}\\
\displaystyle{\| \Lb_b^{\frac{k}{2}}U\|^2_{L^{2;b}}}&\text{if $k$ is even}
\end{cases} \\ 
\nonumber\\
\|U \|^2_{H^{k;b}_{\mathcal L,\me}}&\ge& \lambda_1^{k-j}\|U \|^2_{H^{j;b}_{\mathcal L,\me}}~\ge~ \lambda_1^k\|U\|^2_{L^{2;b}}\quad
\text{for any $j=1,\dots, k$;}
\label{eq:U_embedding}
\end{eqnarray}
\item[$ii)$] the Dirac delta-type function
$$
\delta_0:H^{k;b}_{\mathcal L,\me}(\R\to\mathcal H) \to \mathcal H^{k-\frac{1+b}{2}}_\mathcal L~,\qquad \delta_0(V)=V(0)
$$
is well defined and continuous.
\end{itemize}
\end{Lemma}

\proof
To prove $i)$ use Lemma \ref{L:all}.  Next, let 
$U=\sum\limits_{j=1}^\infty U_j\f_j$ be any curve in $H^{k;b}_{\mathcal L,\me}(\R\to\mathcal L)$. Thanks to $v)$ in Lemma \ref{L:all} we can estimate
$$
\|U\|^2_{H^{k;b}_{\mathcal L,\me}} 
=\sum_{j=1}^\infty~\|U_j\|_{\lambda_j,H^{k;b}_\me}^2\ge m_b \sum_{j=1}^\infty~\lambda_j^{k-\frac{1+b}{2}} |U_j(0)|^2= m_b\|U(0)\|^2_{\mathcal H^{k-\frac{1+b}{2}}_\mathcal L},
$$
which concludes the proof.
\QED

\begin{Remark}
\label{R:U_continuous}
It turns out that $H^{k;b}_{\mathcal L,\me}(\R\to\mathcal H)\subset {\mathcal C}^{0,\frac{1-b_+}{2}}_{\rm loc}(\R\to\mathcal H)$. For the proof, take ${U=}\sum\limits_{j=1}^\infty U_j\f_j  \in H^{1;b}_{\mathcal L,\me}(\R\to\mathcal H)$ and $y_1,y_2\in \R$. We use  
(\ref{eq:loc_Holder}) with $\psi= U_j\in H^{k;b}_\me(\R)$ to estimate
$$
\begin{aligned}
\|U(y_2)-U(y_1)\|_{\mathcal H}^2&=\sum_{j=1}^\infty|U_j(y_2)-U_j(y_1)|^2
&\le  \frac{1}{1-b}\|U\|^2_{H^{1;b}_{\mathcal L,\me}}~\!\big| |y_2|^{1-{b}}-|y_1|^{1-{b}}\big|.
\end{aligned}
$$
Since $H^{k;b}_{\mathcal L,\me}(\R\to\mathcal H)$ is continuously embedded in 
$H^{1;b}_{\mathcal L,\me}(\R\to\mathcal H)$ by (\ref{eq:U_embedding}),  the claim follows.
\end{Remark}

\paragraph{Proof of Theorem \ref{T:draft}}
Recall that ${\mathfrak b}=1-2(s-\lfloor{s}\rfloor)\in(-1,1)$. 
For $u=\sum_ju_j\f_j\in \mathcal H$, we use the notation introduced in (\ref{eq:psilambda})  to rewrite (\ref{eq:best}) as
\begin{equation}
\label{eq:P_simple}
\mathcal P_{\!s}[u](y)= \sum_{j=1}^\infty  \psi_{s,\lambda_j}(y)~\!u_j\f_j~\!.
\end{equation}

We fix any  $u\in \mathcal H^s_\mathcal L$.  
Theorem \ref{T:Aexistence} gives 
$\psi_{s,\lambda_j}\in H^{\lceil{s}\rceil;{\mathfrak b}}_\me(\R)$ and  
$\|\psi_{s,\lambda_j}\|^2_{\lambda_j,H^{\lceil{s}\rceil;{\mathfrak 
b}}_\me}=2d_s\lambda_j^s$ by (\ref{eq:minimum}). 
Thus 
$$
\|\mathcal P_{\!s}[u]\|^2_{H^{\lceil{s}\rceil;{\mathfrak b}}_{\mathcal 
L,\me}}= 
\sum_{j=1}^\infty u_j^2  \|~\! 
\psi_{s,\lambda_j}\|^2_{\lambda_j,H^{\lceil{s}\rceil;{\mathfrak b}}_\me}=
2d_s\sum_{j=1}^\infty \lambda_j^su_j^2 =2d_s\|{\mathcal L}^{\frac{s}{2}}u\|^2_{\mathcal H}=2d_s\|u\|^2_{\mathcal H^s_{\mathcal L}}
$$
 and (\ref{eq:isometry}) is proved. 

\medskip
Next, take any $V\in H^{\lceil{s}\rceil;{\mathfrak b}}_{\mathcal 
L,\me}(\R\to\mathcal H)$ and put $V_j(y)=(V(y),\f_j)_\mathcal H$.
We have
$$
(\mathcal P_{\!s}[u],V)_{\!H^{\lceil{s}\rceil;{\mathfrak b}}_{\mathcal 
L,\me}}=
\sum\limits_{j=1}^\infty u_j(\psi_{s,\lambda_j},V_j)_{\lambda_j, 
H^{\lceil{s}\rceil;{\mathfrak b}}_\me}= 
2d_s \sum\limits_{j=1}^\infty \lambda_j^s u_j V_j(0)=2d_s\langle\mathcal L^su, V(0)\rangle
$$
by (\ref{eq:equation_psi}), which proves $iii)$.  

Evidently $iii)$ implies that $\mathcal P_{\!s}[u]$ is a weak solution to (\ref{eq:ODE}).
Since $\mathcal P_{\!s}[u]$ is smooth on $\R_+$ by Lemma \ref{L:P_smooth}, we see that in fact $\mathcal P_{\!s}[u]$ 
solves (\ref{eq:ODE}) pointwise. The first equality in (\ref{eq:bdry1}) is satisfied by $iii)$ in Lemma \ref{L:P_smooth}.

\medskip 
To conclude the proof of (\ref{eq:bdry1}), we first compute
$$
\Lb_{\mathfrak b}^\is\mathcal P_{\!s}[u](y)=
\sum_{j=1}^\infty \lambda_j^{\is} \big((\Db_{\mathfrak 
b}+1)^{\is}\psi_s\big)(\sqrt{\lambda_j}y)~\!u_j \f_j~\!.
$$
Now we use two items in Lemma \ref{L:old}, namely, $iv)$ (with $m=\is$) and then $i)$ (with $s-\is$ instead of $s$).
This gives
\begin{equation}
\label{eq:brackets}
\begin{aligned}
y^{\mathfrak b}(\partial_y\Lb_{\mathfrak b}^\is\mathcal P_{\!s}[u])(y)&
=\frac{d_s}{d_{s-\is}} 
\sum_{j=1}^\infty \lambda_j^{\is} y^{\mathfrak b} 
\big(\partial_y\psi_{s-\is}\big)(\sqrt{\lambda_j}y)~\!u_j \f_j\\&
= -{d_s} \sum_{j=1}^\infty \psi_{\lceil{s}\rceil-s}(\sqrt{\lambda_j}y)~\!\lambda_j^su_j \f_j=
-d_s\mathcal P_{\!s-\is}[\mathcal L^su](y).
\end{aligned}
\end{equation}
The second limit in (\ref{eq:bdry1}) follows from $iii)$ in Lemma \ref{L:P_smooth},
and $iv)$ is proved.

\medskip

It remains to prove $ii)$. Let  $V\in H^{\lceil{s}\rceil;{\mathfrak 
b}}_{\mathcal L,\me}(\R\to\mathcal H)$ be such that $V(0)=u$. Then
$V_j(0)=u_j$ for any $j\ge 1$. Thus (\ref{eq:minimum}) gives 
$$
u_j^2 \|\psi_{s,\lambda_j}\|^2_{\lambda_j,H^{\lceil{s}\rceil;{\mathfrak 
b}}_\me}
\le \|V_j\|^2_{\lambda_j,H^{\lceil{s}\rceil;{\mathfrak b}}_\me}
$$
for any $j\ge 1$. Thus
$$
\|\mathcal P_{\!s}[u]\|^2_{H^{\lceil{s}\rceil;{\mathfrak b}}_{\mathcal 
L,\me}}=
\sum_{j=1}^\infty u_j^2  \|~\! 
\psi_{s,\lambda_j}\|^2_{\lambda_j,H^{\lceil{s}\rceil;{\mathfrak b}}_\me}
\le \sum_{j=1}^\infty \|V_j\|^2_{\lambda_j,H^{\lceil{s}\rceil;{\mathfrak 
b}}_\me}=
\|V\|^2_{H^{\lceil{s}\rceil;{\mathfrak b}}_{\mathcal L,\me}}.
$$
and $ii)$ follows. The theorem is completely proved.
\QED

\paragraph{Proof of Theorem \ref{T:negative}}
Recall that $\mathcal P_{\!s}[u]: \mathcal H^s_{\mathcal L}\to 
H^{\lceil{s}\rceil;{\mathfrak b}}_{\mathcal L,\me}(\R\to\mathcal H)$
is, up to the constant $2d_s$, an isometry by item $i)$ in Theorem \ref{T:draft}, and that 
$\mathcal L^{-s}: \mathcal H^{-s}_{\mathcal L}\to \mathcal H^s_{\mathcal L}$ is an isometry. 
Thus for any $\zeta\in \mathcal H^{-s}_{\mathcal L}$ we have that
$$
\|\mathcal P_{\!-s}[\zeta]\|_{\!H^{\lceil{s}\rceil;{\mathfrak b}}_{\mathcal 
L,\me}}^2=
2d_s\|\mathcal L^{-s}\zeta\|_{\mathcal H^{s}_{\mathcal L}}= 2d_s\|\zeta\|_{\mathcal H^{-s}_{\mathcal L}},
$$
and (\ref{eq:isometry_negative}) is proved. The conclusions in $iii)$, $iv)$ are immediate consequences of Theorem \ref{T:draft}
(with $u:=\mathcal L^{-s}\zeta$). 

Finally, notice that the strictly convex minimization problem in (\ref{eq:U_minimum_negative}) has a unique solution
$\widehat U\in H^{\lceil{s}\rceil;{\mathfrak b}}_{\mathcal 
L,\me}(\R\to\mathcal H)$, and that $\widehat U$ satisfies
$$
(\widehat U, V)_{\!H^{\lceil{s}\rceil;{\mathfrak b}}_{\mathcal 
L,\me}}=2d_s\langle\zeta, V(0)\rangle= 2d_s\langle\mathcal L^su, 
V(0)\rangle\quad
\text{for any}\quad V\in H^{\lceil{s}\rceil;{\mathfrak b}}_{\mathcal 
L,\me}(\R\to\mathcal H).
$$
Thus $\widehat U= \mathcal P_{\!s}[u]=\mathcal P_{\!-s}[\zeta]$ by $iii)$ in Theorem \ref{T:draft}.
\QED

\section{Generalizations and examples}
\label{S:generalization}

In this section we provide some possible generalizations of our main result. They are based on 
Theorem \ref{T:Aexistence}.

\subsection{Nonnegative operators} 
\label{SS:non-negative}

Assume that $\mathcal L$  is self-adjoint, with a discrete spectrum, nonnegative and with a 
nontrivial kernel. Trivially, for any $s>0$ we have $\ker \mathcal L^s=\ker\mathcal L$, hence
$$
\mathcal L^su =\mathcal L^s(u-\Pi u),
$$
where  $\Pi:\mathcal H\to \mathcal \ker \mathcal L$ is the orthogonal projection on $\ker \mathcal L$.
The domain of the quadratic form $u\mapsto (\mathcal L^s u,u)_\mathcal H$ is
$$
H^s_\mathcal L=\ker\mathcal L\oplus H^s_{\mathcal L_\perp},\quad 
\mathcal L_\perp:= \mathcal L|_{(\ker \mathcal L)^\perp}:(\ker \mathcal L)^\perp\to (\ker \mathcal L)^\perp.
$$
Notice that $\mathcal L_\perp$  is self-adjoint, with a discrete spectrum and positive. Thus 
Theorem \ref{T:draft} provides a full characterization of $\mathcal L^s_\perp$ and of the corresponding quadratic form on $\mathcal H^s_{\mathcal L_\perp}$. This gives, in turn, corresponding results for $\mathcal L^s$ and for its   quadratic form on $\mathcal H^s_{\mathcal L}$.

In particular, the operator $u\mapsto \mathcal P_{\!s}[u]$ in (\ref{eq:best}) is the identity on $\ker\mathcal L$ and 
\begin{equation}
\label{eq:best_kernel}
\mathcal P_{\!s}[u](y)=\Pi[u]+ \mathcal P_{\!s}^\perp[u-\Pi u](y),
\end{equation}
where $\mathcal P^\perp_{\!s}$ is the isometry given by Theorem \ref{T:draft} for the operator   $\mathcal L_\perp$.
Since $\mathcal P_{\!s}[u]$ differs from $\mathcal P_{\!s}^\perp[u-\Pi u]$ by a constant curve, then 
$\mathcal P_{\!s}[u]$, $\mathcal P_{\!s}^\perp[u]$ enjoy the same regularity properties in the Appendix.

\subsection{Non-discrete spectrum}
\label{SS:continuous}

Let $\mathcal L$ be a nonnegative, self-adjoint operator in the Hilbert space $\mathcal H$. Then 
there exists a unique projector-valued spectral measure $E$ on $\R$ supported on the spectrum $\sigma(\mathcal L)\subset [0,\infty)$,
such that
$$
\mathcal L=\int\limits_{[\Lambda,\infty)} \lambda ~\!dE(\lambda),
$$
 where $\Lambda\ge0$ is the bottom of $\sigma(\mathcal L)$ (see e.g., \cite[Ch. 6]{BS}). 
 
For $s>0$, the $s$-power of $\mathcal L$ is formally defined via
$$
\mathcal L^s=\int\limits_{[\Lambda,\infty)} \lambda^s ~\!dE(\lambda).
$$
We denote by $\mathcal H^s_\mathcal L$ the domain of the corresponding quadratic 
form, which is a Hilbert space with norm 
$\|\cdot\|_{\mathcal H^s_\mathcal L}^2= \|\mathcal L^\frac{s}{2}\cdot\|_\mathcal H^2+ \|\cdot\|_\mathcal H^2$.

Let us first assume that $\mathcal L$ be positive definite, i.e. $\Lambda>0$. Then 
$\|\mathcal L^\frac{s}{2}\cdot\|_\mathcal H$ is an equivalent norm in $\mathcal 
H^s_\mathcal L$.

 For $s>0$ non-integer and $ u\in \mathcal H$ we  consider the curve
 \begin{equation}\label{decomp}
\mathcal P_{\!s}[u](y)=\int\limits_{[\Lambda,\infty)}\!\psi_s(\sqrt{\lambda}y)\,
dE(\lambda)u,
\end{equation}
where $\psi_s$ is the function in (\ref{eq:psi_notation}). 
 As in the discrete case, we have that 
$\mathcal P_{\!s}$ maps  any $u\in \mathcal H$ into an even curve in $\mathcal H$; in addition
$\mathcal P_{\!s}[u]\in {\mathcal C}^\infty(\R_+\to \mathcal H^\sigma_\mathcal L)$ for every $u\in \mathcal H, \sigma>0$.

 Further, for $b\in(-1,1)$ we introduce the following (unbounded) operators 
acting on even curves 
 $U\in L^{2;{b}}_\me(\R\to\mathcal H)$,
 $$
\Lb_{ b}U=\int\limits_{[\Lambda,\infty)} (\Db_{ b}+\lambda) 
dE(\lambda)U~\!,\qquad \Db_{ b} U=-\partial^2_{yy}U-{ b}y^{-1}\partial_yU,
$$
compare with (\ref{eq:Db}).

For any integer $k\ge 1$ we introduce the space
$$
H^{k;b}_{\mathcal L,\me}(\R\to\mathcal H)=\Big\{U\in L^{2;b}_\me(\R\to\mathcal H)~\big|~\|  U\|_{H^{k;b}_\me}<\infty~ 
\!~\Big\}. 
$$
Here $\|  \cdot\|_{H^{k;b}_\me}$ is defined 
similarly as we did in the discrete case. 
More precisely, if $k$ is even then
$$
\|  U\|^2_{H^{k;b}_{\mathcal L,\me}}
:=\int\limits_{\R} |y|^b
\Big[\int\limits_{[\Lambda,\infty)} \!
d\big(E(\lambda)V(y,\lambda),V(y,\lambda)\big)\Big]\,dy,  
$$
where $V(y,\lambda)=(\Db_b+\lambda)^{\frac{k}{2}}U(y)$. 
If $k$ is odd then
$$
\|  U\|^2_{H^{k;b}_{\mathcal L,\me}} 
:=
\int\limits_{\R} |y|^b\Big[
\int\limits_{[\Lambda,\infty)} \!
d\big(E(\lambda)\partial_yV(y,\lambda),\partial_yV(y, \lambda)\big)+
\int\limits_{[\Lambda,\infty)} \! \lambda \,
d\big(E(\lambda)V(y,\lambda),V(y,\lambda)\big)\Big]\,dy,
$$
where $V(y,\lambda)=(\Db_b+\lambda)^{\frac{k-1}{2}}U(y)$.

With the above definitions, Theorem \ref{T:draft} holds true, and its proof can be carried out with no essential modifications.
\medskip

If $\Lambda=0$ is an eigenvalue of $\mathcal L$ one can use a decomposition similar to (\ref{eq:best_kernel}) and the above remarks
in the present subsection
for the restriction of $\mathcal L $ to $\ker \mathcal L^\perp$.

A more complicated case is when $0\in \sigma(\mathcal L)$ is not 
an eigenvalue but a point of continuous spectrum. Clearly $\|\mathcal L^\frac{s}{2}\cdot\|_\mathcal H$ cannot bound $\|\cdot\|_\mathcal H$ and therefore 
it is only a seminorm in ${\mathcal H}^s_\mathcal L$. Denote by $\widehat{\mathcal H}^s_\mathcal L$ the completion of 
$\mathcal H^s_\mathcal L$ with respect to $\|\mathcal L^\frac{s}{2}\cdot\|_\mathcal H$.

To avoid additional difficulties, we assume that $\|\mathcal L^\frac{s}{2}\cdot\|_\mathcal H$ is a norm in 
$\widehat{\mathcal H}^s_\mathcal L$.
In this case one can  define a suitable space of curves, and prove a result similar to
Theorem \ref{T:draft}.

\subsection{Examples}
\label{SS:example}

The approach proposed in the present paper can be used, for instance, to recover non-integer powers of 
a large class of differential operators. 

\medskip

The case of the Dirichlet Laplacian in a bounded, smooth domain $\Omega\subset \R^n$ is included in Theorem \ref{T:draft}.
Any curve $y\mapsto U(y)\in L^2(\Omega)=\mathcal H$ is identified with the function $(x,y)\mapsto U(y)(x)$, $\Omega\times\R\to\R$,
so that 
$L^{2;b}(\R\to L^2(\Omega))\equiv L^{2}(\Omega\times\R;|y|^bdxdy)$, and
$$
\|U\|_{L^{2;b}(\R\to L^2(\Omega))}^2=\intr|y|^b\|U(y)\|_{L^2(\Omega)}^2~\!dy=
\iint\limits_{\Omega\times\R}|y|^b|U(x,y)|^2~\!dxdy~\!.
$$
Further, $L^{2;b}_\me(\R\to L^2(\Omega))$ is identified with the space of functions in $L^{2}(\Omega\times\R;|y|^bdxdy)$
which are even in the $y$-variable, that is denoted by $L^{2}_\me(\Omega\times\R;|y|^bdxdy)$. 

We choose $\mathcal L= -\Delta_D$, the standard Laplace operator with domain $H^1_0(\Omega)\cap H^2(\Omega)$. Its 
eigenvalues $\lambda_j$ and corresponding eigenfunctions $\f_j$ solve the Dirichlet problem
$$
\begin{cases}
-\Delta\f_j=\lambda_j\f_j&\text{in $\Omega$}\\
\f_j=0&\text{on $\partial\Omega$,}
\end{cases}~
\qquad \int\limits_\Omega\f_j\f_h~\!dx=\delta_{jh}.
$$
The natural domain $\mathcal H^s_{ -\Delta_D}(\Omega)$ of the quadratic form $u\mapsto ((-\Delta_D)^su,u)_{L^2}$
can be described by the results in \cite[Section 1]{Tr}, see also \cite[Lemma 3]{MN_Mit}:
$$
\mathcal H^s_{ -\Delta_D}(\Omega)=
\displaystyle{\left\{u\in H^s(\Omega)~\Big|~
(-\Delta)^{m} u\big|_{\partial\Omega}=0\ \ 
\text{if}\ \ {m}\in\mathbb N_0, ~ 2{m}<{s-\frac12}
\right\}}
$$
(recall that functions in $H^s(\Omega)$ have a trace on $\partial \Omega$ if and only if $s>\frac12$).

We see that
\begin{equation}
\label{eq:better}
\Lb_{\mathfrak b}U=-{\bf \Delta} U-{\mathfrak b}y^{-1}\partial_yU=-|y|^{-{\mathfrak b}}\div(|y|^{\mathfrak b}\nabla U),
\end{equation}
where  $-{\bf \Delta}$ is the Dirichlet Laplacian in $\Omega\times\R$.

For $s$ non-integer,  Theorem \ref{T:draft} relates the nonlocal operator  $(-\Delta_D)^s$, with the local operator
$\Lb_{\mathfrak b}^{\lceil{s}\rceil}$ acting on  the space $H^{\lceil{s}\rceil;{\mathfrak b}}_{-\Delta_D,\me}(\R\to 
L^2(\Omega))\equiv H^{\lceil{s}\rceil;{\mathfrak b}}_{-\Delta_D,\me}(\Omega\times\R)$. 
For instance, with obvious notation, we have 
$$
\begin{aligned}
H^{1;{\mathfrak b}}_{-\Delta_D,\me}(\Omega\times\R)&= \big\{U\in 
H^{1}_\me(\Omega\times \R;|y|^{\mathfrak b} dxdy)~|~ U(\cdot,y)
\in H^1_0(\Omega)~~\text{for $y\neq 0$}~\big\},\\
\|U\|^2_{H^{1;{\mathfrak b}}_{-\Delta_D,\me}}&=
\iint\limits_{\Omega\times\R}|y|^{\mathfrak b}|\nabla U|^2~\!dxdy~\!;
\\
H^{2;{\mathfrak b}}_{-\Delta_D,\me}(\Omega\times\R)&= 
\big\{U\in H^{1;{\mathfrak b}}_{-\Delta_D,\me}(\Omega\times\R)~|~ 
|y|^{\mathfrak b}\nabla U\in H^{1}(\Omega\times \R;|y|^{-\mathfrak b}dxdy)~\big\},
\\
\|U\|^2_{H^{2;{\mathfrak b}}_{-\Delta_D,\me}}&=
\iint\limits_{\Omega\times\R}|y|^{\mathfrak b}|\bL_{\mathfrak b} U|^2~\!dxdy=
\iint\limits_{\Omega\times\R}|y|^{-\mathfrak b}|\div(|y|^{\mathfrak b}\nabla U)|^2~\!dxdy.
\end{aligned}
$$

\medskip

The Neumann Laplacian in $\Omega$ fits in the situation described in Subsection \ref{SS:non-negative}.
Now we choose $\mathcal L= -\Delta_N$. It is an unbounded operator
on $L^2(\Omega)$ with eigenvalues $\lambda_j\ge 0$ and eigenfunctions $\f_j$ solving
$$
\begin{cases}
-\Delta\f_j=\lambda_j\f_j&\text{in $\Omega$}\\
\partial_\nu\f_j=0&\text{on $\partial\Omega$,}
\end{cases}~
\qquad \int\limits_\Omega\f_j\f_h~\!dx=\delta_{jh}.
$$
The natural domain $\mathcal H^s_{ -\Delta_N}(\Omega)$ of the quadratic form $u\mapsto ((-\Delta_N)^su,u)_{L^2}$
is
$$
\mathcal H^s_{ -\Delta_N}(\Omega)=
\displaystyle{\left\{u\in H^s(\Omega)~\Big|~
\partial_\nu(-\Delta)^{m} u\big|_{\partial\Omega}=0\ \ 
\text{if}\ \ {m}\in\mathbb N_0, ~ 2{m}<{s-\frac32}
\right\}},
$$
see \cite[Section 1]{Tr}.

In this case, the operator $\Lb_{\mathfrak b}$ is pointwise defined as in (\ref{eq:better}).
For $s\notin \mathbb N$, the nonlocal operator  $(-\Delta_N)^s$ is related to   $\Lb_{\mathfrak b}^{\lceil{s}\rceil}$, 
acting on a different domain 
$$H^{\lceil{s}\rceil;{\mathfrak b}}_{-\Delta_N,\me}(\R\to L^2(\Omega))\equiv H^{\lceil{s}\rceil;{\mathfrak b}}_{-\Delta_N,\me}(\Omega\times\R).
$$
The approach described in  Subsection \ref{SS:non-negative} covers this example, as $\lambda_1=0$.

\medskip

Lastly, if $n>2s$ then the fractional Laplacian $\Ds$ on $\R^n$ fits into the general approach in Subsection \ref{SS:continuous}. In this case, thanks to Hardy inequality the space $\widehat{\mathcal H}^s_{-\Delta}$
can be identified with the standard homogeneous Sobolev space 
$\mathcal D^s(\R^n)\hookrightarrow L^2(\R^n;|x|^{-2s}dx)$. The resulting space of curves can be identified with the space
$\mathcal D^{\lceil{s}\rceil; \mathfrak b}_\me(\R^{n+1})$ in \cite{CM}.


\appendix
\section{On the transforms $\mathcal P_{\!s}$} 
\label{A}

Here we assume that $s>0$ is non-integer and study the transform $\mathcal P_{\!s}[\cdot]$, see (\ref{eq:P_simple}).
We start by noticing that formulae (\ref{eq:dual0}) and (\ref{eq:Hs_norm}) hold and that 
$\mathcal H^s_{\mathcal L}$ is the domain of the quadratic form of $\mathcal L^s$,
for negative orders $s$ as well.

\begin{Lemma}
\label{L:P_smooth}
Let $s>0$, $\sigma\in\R$.  
\begin{itemize}
\item[$i)$] For any $u\in\mathcal H$, we have $\mathcal P_{\!s}[u]\in {\mathcal C}^\infty(\R_+\to \mathcal H^\sigma_\mathcal L)$, and
$\|\partial^k_y\mathcal P_{\!s}[u](y)\|_{\mathcal H^\sigma_\mathcal L}$ 
decays exponentially 
 as $y\to\infty$, for any order $k\ge 0$;
 \item[$ii)$] The linear operator $u\mapsto \mathcal P_{\!s}[u](y)$ is 
nonexpansive in $H^\sigma_\mathcal L$, that is,
\begin{equation}
\label{eq:nonexpansiveA}
\|\mathcal P_{\!s}[u](y)\|_{\mathcal H^{\sigma}_\mathcal L}\le \|u\|_{\mathcal H^\sigma_\mathcal L}
\quad \text{for any $y\in\R$;}
\end{equation}
 \item[$iii)$] If $u\in \mathcal H^\sigma_\mathcal L$ then
  $\mathcal P_{\!s}[u]\in {\mathcal C}^0(\R\to \mathcal H^\sigma_\mathcal L)$ and 
 $\mathcal P_{\!s}[u](0)=u$;
\item[$iv)$]
The  operator $u\mapsto \mathcal P_{\!s}[u](y)$ commutes with 
the fractional powers of $\mathcal L$, that is,
\begin{gather}
\label{eq:commute}
\mathcal P_{\!s}[\mathcal L^\sigma u](y)= \mathcal L^\sigma \big(\mathcal P_{\!s}[u](y)).
\end{gather}
\end{itemize}
\end{Lemma}

\proof 
By the properties of the Bessel functions, for any 
integer $k\ge 0$ and any $\delta>0$ we have $|(\partial^k_y\psi_s)(y)|\le c(\delta)e^{-y}$ for $y>\sqrt{\lambda_1}\delta$, where the constant $c(\delta)$
depends on $\delta, s$ and $k$ but not on $y$. Thus, for $y>\delta$ we have
$$
\lambda_j^{k+\sigma} |(\partial^k_y\psi_s)(\sqrt{\lambda_j}y)|^2\le c(\delta)^2 \lambda_j^{k+\sigma} e^{-2\sqrt{\lambda_j}y}\le C(\delta)e^{-\sqrt{\lambda_1}y}
$$
because $\lambda_j\ge \lambda_1>0$, where the new constant $C(\delta)$ depends only on $\delta, s, \sigma$ and $ k$. 
It readily follows that 
$$
\|\partial^k_y\mathcal P_{\!s}[u](y)\|^2_{\mathcal H^\sigma_\mathcal L}=\sum_{j=1}^\infty \lambda_j^{k+\sigma} u_j^2 |(\partial^k_y\psi_s)(\sqrt{\lambda_j}y)|^2
\le C(\delta) \|u\|^2_{\mathcal H} ~\!e^{-\sqrt{\lambda_1}y}
$$
for any $u\in\mathcal H$, provided that $y>\delta$, and $i)$ is proved.  

\medskip 
Now we take $u\in \mathcal H^\sigma_\mathcal L$. 
By Remark \ref{R:max_psi}, we have $0< \psi_{s,\lambda_j}(y)\le \psi_{s,\lambda_j}(0)=1$. Thus
$$
\|\mathcal P_{\!s}[u](y)\|^2_{\mathcal H^{\sigma}_\mathcal L}=
\sum_{j=1}^\infty  \lambda_j^\sigma u_j^2 (\psi_{s,\lambda_j}(y))^2 \le
\sum_{j=1}^\infty  \lambda_j^\sigma u_j^2 =\|u\|^2_{\mathcal H^\sigma_\mathcal L}~\!,
$$
which proves $ii)$. Further, we have
\begin{equation}
\label{eq:lebesgue0}
\|u-\mathcal P_{\!s}[u](y)\|^2_{\mathcal H^\sigma_\mathcal L}=\sum_{j=1}^\infty \lambda^\sigma_j u_j^2(\psi_{s,\lambda_j}(0)-\psi_{s,\lambda_j}(y))^2
\le \sum_{j=1}^\infty \lambda^\sigma_ju_j^2.
\end{equation}
The first series in (\ref{eq:lebesgue0}) is dominated by a convergent number series and converges to zero termwise
as $y\to 0$. We infer that $\|u-\mathcal P_{\!s}[u](y)\|^2_{\mathcal H^\sigma_\mathcal L}\to 0$ as $y\to 0$, which implies
$\mathcal P_{\!s}[u]\in {\mathcal C}^0(\R\to \mathcal H^\sigma_\mathcal L)$, and $iii)$ is proved.

Since the equality in $iv)$ is trivial, the proof is complete.
\QED

Thanks to Lemma \ref{L:P_smooth} and (\ref{eq:brackets}), we can improve the convergences in (\ref{eq:bdry1}) as follows.

\begin{Corollary}
\label{C:convergence}
Let $s>0$ be non-integer, ${\mathfrak b}=1-2(s-\is)$. Assume that 
$u\in \mathcal H^\sigma_\mathcal L$ for some $\sigma\in\R$.
Then $\mathcal P_{\!s}[u]$ solves the differential equation (\ref{eq:ODE}) and satisfies the boundary conditions
$$
\begin{aligned}
\lim\limits_{y \to 0^+}\mathcal P_{\!s}[u](0)&= u\quad\text{in~ $\mathcal H^{\sigma}_{\mathcal L}$}\\
\lim\limits_{y \to 0^+}y^{\mathfrak b}~\!{\partial_y} 
\big(\bL_{\mathfrak b}^{\lfloor{s}\rfloor}\mathcal P_{\!s}[u]\big)(y)&= - 
d_s\,{\mathcal L}^{s} u
\quad\text{in~ $\mathcal H^{\sigma-2s}_{\mathcal L}$.}
\end{aligned}
$$
\end{Corollary}

\begin{Remark}
\label{R:half_integer}
For any integer $k\ge 0$ we have
$$
\psi_{k+\frac12}(y)=\frac{1}{(k+1)!}|y|^{k}e^{-|y|},\qquad \mathcal P_{\!k+\frac12}[u](y)=\frac{1}{(k+1)!}|y|^k~\!\mathcal P_{\!\frac12}[\mathcal L^\frac12 u](y)~\!.
$$
\end{Remark}

\subsection{Derivatives}

The regularity of the curve $\mathcal P_{\!s}[u]$ given  in Lemma \ref{L:P_smooth}  
improves as $s$ increases.  
We start by proving a technical result which involves the Beta function
 $$
 \textrm{B}(\tau,t)=\int\limits_0^1x^{\tau-1} (1-x)^{t-1}~\!dx=\frac{\Gamma(t)\Gamma(\tau)}{\Gamma(t+\tau)}~\!.
 $$
The coefficients in the next lemma are computed by taking inspiration from \cite[Section~4]{CM}.

\begin{Lemma}
\label{L:3}
Let $\sigma\in\R$,  $u\in\mathcal H^\sigma_\mathcal L$, $y>0$.
\begin{itemize}
\item[$i)$] If $s\in(0,1)$ then 
$\partial_{y}\mathcal P_{\!s}[u](y)=
-d_s y^{2s-1} \mathcal P_{\!1-s}[\mathcal L^su](y)$;
\item[$ii)$] If $s>1$ then for any  $m=1,\dots, \is$ it holds that 
\end{itemize}
\begin{gather}
\label{eq:even_d}
\partial^{2m}_{y}\mathcal P_{\!s}[u](y)=\frac{1}{{\textrm{\rm B}\big(s,\frac12\big)}}
\sum_{\ell=0}^{m} \binom{{m}}{\ell} {(-1)^{\ell} 
\textrm{\rm B}\big(s-\ell,\tfrac12\big)}\cdot\mathcal P_{{s}-\ell}[\mathcal L^{m}u](y)\\
\label{eq:odd_d}
\partial^{2{m}-1}_y \mathcal P_{{s}}[u](y) =
y\cdot \frac{1}{{\textrm{\rm B}\big(s,\frac12\big)}}
\sum_{\ell=1}^{m} \binom{{m-1}}{\ell-1} (-1)^{\ell} 
{\textrm{\rm B}\big(s-\ell,\tfrac32\big)}\cdot\mathcal P_{{s}-\ell}[\mathcal L^{m}u](y).
\end{gather}
\end{Lemma}

\proof
If $s\in(0,1)$ then $\partial_y\psi_{s,\lambda}(y)=-d_s y^{2s-1} \lambda^s\psi_{1-s,\lambda}(y)$ by Lemma 
\ref{L:old}. Thus
$$
\partial_{y}\mathcal P_{\!s}[u](y)=
-d_sy^{2s-1}\cdot \sum_{j=1}^\infty \partial^{2m}_{y}\psi_{1-s,\lambda_j}(y)~\!\lambda_j^su_j \f_j,
$$
and the identity in $i)$ follows.

\medskip
To handle the case $s>1$ we put
$\displaystyle{
\gamma_{s,\ell}= \frac{\textrm{\rm B}(s-\ell,\frac12)}{\textrm{\rm B}(s,\frac12)}}
=\frac{\Gamma\big({s}+\frac12\big)}{\Gamma({s})}\frac{\Gamma({s}-\ell)}{\Gamma\big({s}+\frac12-\ell\big)}$.
Using $ii)$ in Lemma \ref{L:old}
and induction one gets
\begin{equation}
\label{eq:even_psi}
\partial^{2m}_{y}\psi_s(y)=
\sum_{\ell=0}^{m} \binom{{m}}{\ell} (-1)^{\ell} \gamma_{s,\ell} ~\psi_{{s}-\ell}(y).
\end{equation}
for any integer $m=1,\dots,\is$. Since $\partial^{2m}_{y}\psi_{s,\lambda}(y)=\lambda^m\big(\partial^{2m}_{y}\psi\big)(\sqrt\lambda y)$, we infer that
$$
\partial^{2m}_{y}\mathcal P_{\!s}[u](y)=
\sum_{j=1}^\infty  \partial^{2m}_{y}\psi_{s,\lambda_j}(y)~\!u_j\f_j
= \sum_{\ell=0}^{m} \binom{{m}}{\ell} (-1)^{\ell} \gamma_{s,\ell} 
\sum_{j=1}^\infty \psi_{{s}-\ell,\lambda_j}(y)~\!\lambda_j^m u_j \f_j,
$$
which proves (\ref{eq:even_d}). 

Arguing as for $i)$ we obtain
$$
\partial_{y}\mathcal P_{\!s}[u](y)=
-\frac{y}{2(s-1)} \mathcal P_{\!s-1}[\mathcal Lu](y),
$$
i.e. (\ref{eq:odd_d}) holds if $m=1$.
If $m>1$ we use (\ref{eq:even_psi}) for $m-1$ and then $i)$ in Lemma \ref{L:old} to compute
$$
\begin{aligned}
\partial^{2{m}-1}_y \psi_s(y)&= \partial_y~\!\partial^{2({m}-1)}_y \psi_s= 
\sum_{\ell=0}^{m-1}  \binom{{m-1}}{\ell} (-1)^{\ell}\gamma_{s,\ell} ~\partial_y\psi_{{s}-\ell}(y)\\ 
&
=y\cdot \sum_{\ell=0}^{m-1} \binom{{m-1}}{\ell} (-1)^{\ell+1} \frac{\gamma_{s,\ell}}{2(s-\ell-1)} ~
\psi_{{s}-\ell-1}(y)\\
&=
y\cdot \sum_{\ell=1}^{m}  \binom{{m}}{\ell}(-1)^{\ell} \frac{\ell~\!\gamma_{s,\ell}}{2m(s+\frac12-\ell)}
~\!\psi_{{s}-\ell}(y)~\!.
\end{aligned}
$$
Then (\ref{eq:odd_d}) follows by arguing as in the "even" case.
\QED

\begin{Theorem}
\label{T:bounded_derivatives}
Let $2s\ge 1$, $\sigma\in\R$ and let $k$ be an integer, with $1\le k\le \istwo$.
\begin{itemize}
\item[$i)$] Let $u\in \mathcal H^\sigma_\mathcal L$. Then 
\begin{equation}
\label{eq:derivative_estimate}
\|\partial^k_y\mathcal P_{\!s}[u](y)\|_{\mathcal H^{\sigma-k}_\mathcal L}\le c_k\|u\|_{\mathcal H^\sigma_\mathcal L}
\quad \text{for any $y>0$,}
\end{equation}
where the constant $c_k$ depends only on $s$ and $k$. Thus, for any $y>0$ the linear operator $u\mapsto \partial^k_y\mathcal P_{\!s}[u](y)$ is 
continuous  $\mathcal H^\sigma_\mathcal L\to \mathcal H^{\sigma-k}_\mathcal L$;
 \item[$ii)$] If in addition\footnote{\footnotesize{this is a restriction only if $s$ is a half integer}} $k<2s$ then 
 $\partial^k_y\mathcal P_{\!s}[u]\in {\mathcal C}^0(\R\to \mathcal H^{\sigma-k}_\mathcal L)$ for any
 $u\in \mathcal H^\sigma_\mathcal L$.
  \end{itemize}
\end{Theorem}

\proof
It is convenient to define
$$
M_{\alpha,\beta}=\max_{y\ge 0} y^{2\beta}\psi_\alpha(y)^2~,\quad \alpha, \beta>0.
$$
If $\frac12\le s<1$, then $i)$ in Lemma \ref{L:3} gives
$$
\|\partial_{y}\mathcal P_{\!s}[u](y)\|^2_{\mathcal H^{\sigma-1}_\mathcal L}= 
d_s^2 \|y^{2s-1}\mathcal P_{\!1-s}[\mathcal L^su](y)\|^2_{\mathcal H^{\sigma-1}_\mathcal L}.
$$
The conclusion in $i)$ follows, because
\begin{equation}
\label{eq:name}
\begin{aligned}
 \|y^{2s-1}\mathcal P_{\!1-s}[\mathcal L^su]&(y)\|^2_{\mathcal H^{\sigma-1}_\mathcal L}=
 \sum_{j=1}^\infty \lambda_j^{\sigma-1} y^{2(2s-1)} \lambda_j^{2s} u_j^2 \psi_{1-s}(\sqrt{\lambda_j}y)^2\\ &
=\sum_{j=1}^\infty \lambda_j^{\sigma} u_j^2(\sqrt{\lambda_j}y)^{2(2s-1)}\psi_{1-s}(\sqrt{\lambda_j}y)^2
\le M_{1-s,2s-1}\|u\|^2_{\mathcal H^{\sigma-1}_\mathcal L}.
\end{aligned}
\end{equation}
If $2s>1$, then the series in (\ref{eq:name}) are dominated by a convergent number series 
and converge to zero termwise as $y\to 0$. We infer that $\partial_{y}\mathcal P_{\!s}[u](y)\to 0$
in $\mathcal H^{\sigma-1}_\mathcal L$ as $y\to 0$, which proves $ii)$ in this case.

Next, let $s>1$. We first face the case when $k\le 2\is$ is even. 
Take integers $\ell, m$ with $0\le \ell\le m\le \is$. By Lemma \ref{L:P_smooth}
we have 
$$
\|\mathcal P_{{s}-\ell}[\mathcal L^m u](y)\|_{\mathcal H^{\sigma-2m}_\mathcal L}\le
\|\mathcal L^m u\|_{\mathcal H^{\sigma-2m}_\mathcal L}=
\|u\|_{\mathcal H^{\sigma}_\mathcal L}~,\quad \mathcal P_{{s}-\ell}[\mathcal L^m u]\in \mathcal C^0(\R\to \mathcal H^{\sigma-2m}_\mathcal L).
$$
Taking also (\ref{eq:even_d}) into account, we see that the conclusions hold in this case.

Let now $k\le 2\is-1$ be odd. For $1\le\ell\le m$ 
we estimate
\begin{equation}
\label{eq:name2}
\begin{aligned}
\|y \mathcal P_{{s}-\ell}[\mathcal L^m u](y)\|^2_{\mathcal H^{\sigma-2m+1}_\mathcal L}
&=
\|y \mathcal L^m(\mathcal P_{{s}-\ell}[u](y))\|^2_{\mathcal H^{\sigma-2m+1}_\mathcal L}
=
\|y \mathcal  P_{{s}-\ell}[u](y)\|^2_{\mathcal H^{\sigma+1}_\mathcal L}\\
&=
\sum_{j=1}^\infty \lambda_j^\sigma u_j^2 (\sqrt{\lambda_j}y)^2\psi_{s-\ell}(\sqrt{\lambda_j}y)^2
\le
M_{s-\ell,1}\|u\|^2_{\mathcal H^{\sigma}_\mathcal L}.
\end{aligned}
\end{equation}
In view of (\ref{eq:odd_d}), we see that (\ref{eq:derivative_estimate}) holds also in this case.
By repeating the argument we used for $\frac12\le s<1$ one plainly conclude the proof also in this case.

It remains to discuss the case $\is+\frac12\le s<\lceil{s}\rceil$ and $k=2\is+1=\istwo$.
We differentiate formula (\ref{eq:even_d}) for $m=\is$.
To compute $\partial_y \mathcal P_{{s}-\ell}[\mathcal L^{\is}u](y)$,
 we use (\ref{eq:odd_d})
for $\ell=1,\dots,\is-1$ and $i)$ in Lemma \ref{L:3} for the last $\ell$. It gives
\begin{equation}
\label{eq:last}
\begin{aligned}
\partial^{2{\is}+1}_y \mathcal P_{{s}}[u](y) &=- 
\sum_{\ell=1}^{\is} {a}_{s,\ell} \cdot\big(y \mathcal P_{{s}-\ell}[\mathcal L^{\lceil{s}\rceil}u](y)\big)\\ &
-{a}_{s}\cdot\big(y^{2(s-\is)-1}  ~\mathcal P_{\lceil{s}\rceil-s}[\mathcal L^{s}u](y)\big).
\end{aligned}
\end{equation}
where the coefficients $a_{s,\ell}, a_s\in\R$ depend only on $s$ and $\ell$.
One can easily adapt the arguments we used for 
(\ref{eq:name2}), (\ref{eq:name}). In this way one proves $i)$ if $\is+\frac12\le s<\lceil{s}\rceil$, and $ii)$
if $\is+\frac12< s<\lceil{s}\rceil$.
\QED

\begin{Theorem}
\label{C:3}
Let $s>1$, $\sigma\in\R$, $u\in\mathcal H^\sigma_\mathcal L$.
Then, for any $k=1,\dots,\is$ we have  
\begin{equation}
\label{eq:Taylor}
\mathcal P_{\!s}[u](y)=\frac{1}{\Gamma(s)}
\sum_{m=1}^k \frac{\Gamma(s-m)}{2^{2m}m!}\cdot\mathcal L^m u\cdot y^{2m}+ {o(y^{2k}})\qquad
\text{as $y\to 0$}
\end{equation}
with convergence  in $\mathcal H^{\sigma-2k}_\mathcal L$.
\end{Theorem}

\proof 
Take an integer $k=1,\dots,\is$. 
By $ii)$ in Theorem \ref{T:bounded_derivatives} we have that $\mathcal P_{\!s}[u]\in \mathcal C^{2k}
(\R\to \mathcal H^{\sigma-2k}_\mathcal L)$. Further, for any $m=1,\dots,k$, Lemma \ref{L:3} gives
$$
\partial^{2m}_{y}\mathcal P_{\!s}[u](0)
=
\frac{1}{\textrm{ B}(s,\frac12)}
 \sum\limits_{\ell=0}^{m} \binom{{m}}{\ell}  (-1)^{\ell} 
{\textrm{\rm B}(s-\ell,\tfrac12)}\cdot\mathcal L^{m}u
$$
and $\partial^{2m-1}_{y}\mathcal P_{\!s}[u](0)=0$. 
Then (\ref{eq:Taylor}) follows via  Taylor expansion formula, thanks to Lemma \ref{L:thanks} below.
\QED

\begin{Lemma}
\label{L:thanks}
Let $m\le\is$ be a positive integer. Then
$$
\kappa_{s,m}:= \frac{1}{\textrm{ B}(s,\frac12)}
 \sum\limits_{\ell=0}^{m} \binom{{m}}{\ell}  (-1)^{\ell} 
{\textrm{\rm B}(s-\ell,\tfrac12)}= (-1)^m\frac{\Gamma(s-m)}{\Gamma(s)}\frac{1}{2^{2m} m!}~\!(2m)!~\!.
$$
\end{Lemma}

\proof
We compute
  $$
  \begin{aligned}
  \sum\limits_{\ell=0}^{m}  \binom{{m}}{\ell}  (-1)^{\ell}&
  \textrm{B}\big({s}-\ell,\tfrac12\big)=
  \int\limits_0^1 x^{-\frac12} (1-x)^{s-m-1}
  \Big( \sum\limits_{\ell=0}^{m} \binom{{m}}{\ell}  (-1)^{\ell} (1-x)^{m-\ell}\Big)dx\\
  &=
  (-1)^m \int\limits_0^1 x^{m-\frac12} (1-x)^{s-m-1}dx=
  (-1)^m \textrm{B}(s-m,m+\tfrac12).
  \end{aligned}
  $$
  Recalling the  Legendre duplication formula, we infer that
  $$
\begin{aligned}
\kappa_{s,m}&= (-1)^m \frac{\textrm{B}(s-m,m+\frac12)}{\textrm{B}(s,\frac12)}=
(-1)^m\frac{\Gamma(s-m)}{\Gamma(s)}~\!
\frac{\Gamma(m+\frac12)}{\sqrt{\pi}}\\ &=
(-1)^m\frac{\Gamma(s-m)}{\Gamma(s)}~\!\frac{2^{1-2m}\Gamma(2m)}{\Gamma(m)}
= (-1)^m\frac{\Gamma(s-m)}{\Gamma(s)}\frac{1}{2^{2m} m!}~\!(2m)!,
\end{aligned}
$$
which completes the proof.
 \QED

\begin{Corollary}
\label{C:bdry}
Let $s>1$, $u\in \mathcal H^s_\mathcal L$. Then
for any integer $m=1,\dots,\is$ we have that
$$
\bL_b^{m}\mathcal P_{\!s}[u](y) = \frac{d_s}{d_{s-m}} \mathcal P_{\!s-m}[\mathcal L^mu](y)~,\quad  y\in \R.
$$
\begin{gather*}
{\lim\limits_{y\to 0} y^{-1}{\partial^{2m-1}_y}\mathcal P_{\!s}[u]} =
\lim\limits_{y\to 0}
{{\partial^{2m}_y}\mathcal P_{\!s}[u] }= \kappa_{s,{m}} \mathcal L^{m} u
\end{gather*}
where $\kappa_{s,m}$ is the  constant in Lemma 
\ref{L:thanks}. The  limits are taken in the $\mathcal H^{s-2m}_\mathcal L$ topology. 
\end{Corollary}

\proof
The first equality follows from formulae (\ref{eq:rescaling100}) and (\ref{eq:Abm}):
$$
\begin{aligned}
\bL_b^{m}\mathcal P_{\!s}[u](y)&=\sum_{j=1}^\infty u_j (\Db_b+\lambda_j)^m\psi_{s,{\lambda_j}}(y)~\!u_j\f_j
= \sum_{j=1}^\infty \lambda_j^m [(\Db_b+1)^m\psi_{s}](\sqrt\lambda_jy)~\!u_j\f_j\\
&=\frac{d}{d_{s-m}} \sum_{j=1}^\infty  \psi_{s-m,{\lambda_j}}(y)~\!\lambda_j^m u_j\f_j.
\end{aligned}
$$
To conclude the proof, use $ii)$ in Lemma \ref{L:3} and then $iii)$ in Lemma \ref{L:P_smooth}.
\QED

Our last result in this section involves the  H\"older-type spaces $\widetilde{\mathcal C}^{\alpha}$ in (\ref{eq:Holder_space}).

\begin{Theorem}
\label{T:4}
Let $s>0$ {non-integer}, 
$\sigma\in\R$, $u\in \mathcal H^\sigma_\mathcal L$, $\alpha\in(0,2s]$. Then 
\begin{gather}
\label{eq:T_holder}
\mathcal P_{\!s}[u] \in \widetilde{\mathcal C}^{\alpha}(\R\to \mathcal H^{\sigma-\alpha}_\mathcal L)
~,\qquad \llbracket \mathcal P_{\!s}[u]\rrbracket_{\widetilde{\mathcal C}^{\alpha}}\le c\|u\|_{\mathcal H^\sigma_\mathcal L}.
\end{gather}
\end{Theorem}

\proof
Thanks to $ii)$ in Theorem \ref{T:bounded_derivatives},  we only have to investigate the H\"olderianity of $\partial_y^\isalpha \mathcal P_{\!s}[u]$ if $\alpha>\isalpha$, and the Lipschitz properties
of $\partial_y^{\alpha-1} \mathcal P_{\!s}[u]$ if $\alpha$ is integer.

Theorem \ref{T:Aexistence} already gives $\psi_s\in\widetilde{\mathcal C}^{2s}(\R)$. Since $\psi_s$
decays exponentially at infinity together with its derivatives of any order, we infer that 
$\psi_s\in\widetilde{\mathcal C}^\alpha(\R)$ for any $\alpha\in(0,2s]$. 
Since trivially $\partial^{k}_y\psi_{s,\lambda}(y)=\lambda^\frac{k}2 (\partial^{k}_y\psi_{s})(\sqrt{\lambda}y)$
for any integer $k$ and any $\lambda>0$, then 
$\llbracket \psi_{s,\lambda}\rrbracket_{\widetilde{\mathcal C}^{\alpha}} =  
\lambda^{\frac\alpha2} \llbracket \psi_{s}\rrbracket_{\widetilde{\mathcal C}^{\alpha}}$ for any $\alpha\in(0,2s]$.  

Take arbitrary points $y_1, y_2\in\R$. Without loss of generality, we can assume
that $y_1, y_2\ge 0$. If $\alpha$ is not an integer, then 
$$
\begin{aligned}
\|\partial^{\isalpha}_y\mathcal P_{\!s}&[u](y_1)-\partial^{\isalpha}_y\mathcal P_{\!s}[u](y_2)\|_{\mathcal H^{\sigma-\alpha}_\mathcal L}^2
=
\sum_{j=1}^\infty \lambda_j^{\sigma-\alpha}u_j^2 | \partial^{\isalpha}_y\psi_{s,\lambda_j}(y_1)-\partial^{\isalpha}_y\psi_{s,\lambda_j}(y_2)|^2\\
&\le
\llbracket \psi_{s}\rrbracket_{\widetilde{\mathcal C}^{\alpha}}^2
\sum_{j=1}^\infty \lambda_j^{\sigma-\alpha} u_j^2\lambda^{\alpha}|y_1- y_2|^{2(\alpha-\isalpha)}
=\llbracket \psi_{s}\rrbracket^2_{\widetilde{\mathcal C}^{\alpha}}\|u\|^2_{\mathcal H^\sigma_\mathcal L}|y_1- y_2|^{2(\alpha-\isalpha)}.
\end{aligned}
$$
If $\alpha$ is integer, with a similar computation we get
$$
\begin{aligned}
\|\partial^{\alpha-1}_y\mathcal P_{\!s}[u](y_1)&-\partial^{\alpha-1}_y\mathcal P_{\!s}[u](y_2)\|_{\mathcal H^{\sigma-\alpha}_\mathcal L}^2
\le 
c\sum_{j=1}^\infty \lambda_j^{\sigma} u_j^2|y_1- y_2|^{2}
=c\|u\|^2_{\mathcal H^\sigma_\mathcal L}|y_1- y_2|^2.
\end{aligned}
$$
In both cases, this concludes the proof.
\QED

\subsection{Isometric properties}
\label{SS:isometric}

From Theorem \ref{T:draft} we already know that the linear transform $u\mapsto \mathcal P_{\!s}[u]$ is, up to a constant,
an isometry $\mathcal H^s_\mathcal L \to H^{\lceil{s}\rceil;{\mathfrak 
b}}_{\mathcal L,\me}(\R\to\mathcal H)$ for
${\mathfrak b}:=1-2(s-\is)$.
In this section we point out more isometric properties of $\mathcal P_{\!s}$. We stress the fact that
$s>0$ might be an integer number.

\begin{Theorem}
\label{P:jump}
{Let $s>0$, $b\in(-1,1)$ and $\sigma\in \R$}. Up to a  constant (not 
depending on $\sigma$), the operator 
$\mathcal P_{\!s}$ is an isometry 
$\mathcal H^\sigma_\mathcal L\to L^{2;b}_\me(\R\to \mathcal 
H^{\sigma+\frac{1+b}{2}}_{\mathcal L})$. More precisely,
\begin{equation}
\label{eq:iso1}
\|\mathcal P_{\!s}[u]\|_{L^{2;b}(\R\to \mathcal 
H^{\sigma+\frac{1+b}{2}}_{\mathcal L})}
=\|\psi_s\|_{L^{2;b}(\R)}~\!\|u\|_{\mathcal H^\sigma_\mathcal 
L}\qquad\text{for any $u\in \mathcal H^\sigma_\mathcal L$.}
\end{equation}
\end{Theorem}

\proof For $u\in\mathcal  H^\sigma_\mathcal L$ we compute
$$
\begin{aligned}
\intr|y|^b\| \mathcal P_{\!s}[u](y)&\|_{\mathcal 
H^{\sigma+\frac{1+b}{2}}}^2dy= \sum_{j=1}^\infty 
\lambda_j^{\sigma+\frac{1+b}{2}}  u_j^2~\! 
\intr|y|^b|\psi_{s}(\sqrt{\lambda_j}y)|^2~dy\\
&= \big(\intr|y|^b|\psi_{s}(y)|^2~dy\big)\sum_{j=1}^\infty  
\lambda_j^\sigma u_j^2= 
\big(\intr|y|^b|\psi_{s}(y)|^2~dy\big)\|u\|^2_{\mathcal H^\sigma_\mathcal 
L}~\!,
\end{aligned}
$$
and the Lemma is proved.
\QED

Let $\alpha> 0$. 
We recall the definition of the Sobolev--Slobodetskii spaces and corresponding seminorms
$$
\begin{aligned}
H^\alpha(\R)&=\{ \psi\in L^2(\R)~|~ \llbracket \psi\rrbracket^2_{H^\alpha}<\infty\}\\
\llbracket \psi\rrbracket^2_{H^\alpha}&=\intr |(-\partial^2_{yy})^\frac\alpha2\psi(y)|^2dy
=\intr |\xi|^{2\alpha}|\widehat \psi(\xi)|^2d\xi~\!,
\end{aligned}
$$
where $\widehat{\psi}$ stands for the unitary Fourier transform of $\psi\in L^2(\R)$, namely, 
$$
\widehat{\psi}(\xi)=\frac{1}{\sqrt{2\pi}}\intr e^{-ixy}\psi(y)dy~\!.
$$

We first compute the Fourier transform of the function $\psi_s$ in (\ref{eq:psi_notation}).

\begin{Proposition}
\label{P:psi_transform}
Let $s>0$ (possibly integer). Then
$$
\widehat{\psi_s}(\xi)=\frac{\sqrt2 \Gamma\big(s+\frac12\big)}{\Gamma(s)}~\!(1+\xi^2)^{-\frac{1+2s}{2}}.
$$
In particular, 
$\psi_s\in H^\alpha(\R)$ if and only if $\alpha<2s+\frac12$, and in this case
\begin{equation}
\label{eq:psi_norms}
\llbracket \psi_s\rrbracket^2_{H^\alpha}=
\frac{\Gamma\big(s+\tfrac12 \big)^2}{s\Gamma(2s)\Gamma(s)^2}~\Gamma\big(\alpha+\tfrac12 \big) 
\Gamma\big(2s-\alpha+\tfrac12 \big).
\end{equation}
\end{Proposition}

\proof
It is well known, see for instance \cite[Lemma 4.2]{CM} for a simple proof, that
$$
\reallywidehat{(1+|\cdot|^2)^{-\frac{1+2s}{2}}}(y)= \frac{\Gamma(s)}{\sqrt2 \Gamma\big(s+\frac12\big)}~\!\psi_s(y).
$$
To conclude, use the symmetry of $\psi_s$ and make direct computations.
\QED

For $\alpha>0$ we introduce a Sobolev-type space of curves $\R\to \mathcal H$ and corresponding 
seminorm as follows:
\begin{gather*}
H^\alpha(\R\to \mathcal H)=\big\{ U\in L^2(\R\to\mathcal H)~\big|~ 
\llbracket U\rrbracket^2_{H^\alpha}:=\sum_{j=1}^\infty\intr |\xi|^{2\alpha}|\widehat U_j(\xi)|^2d\xi<\infty\big\},
\end{gather*}
where the Fourier transform of a function $U=\sum_jU_j\f_j\in L^2(\R\to\mathcal H)$
is defined via the Fourier transform of its coordinates, that is,
$$
\widehat U(\xi)=\sum_{j=1}^\infty \widehat U_j(\xi)\f_j.
$$
It is evident that $H^\alpha(\R\to \mathcal H)$ is a Hilbert space with norm 
$$
\|U\|^2_{H^\alpha}=\llbracket U\rrbracket^2_{H^\alpha}+ \|U\|^2_{L^2}
=\sum_{j=1}^\infty \intr(|\xi|^{2\alpha}+1)|\widehat U_j(\xi)|^2~\!d\xi.
$$

\begin{Theorem}
\label{T:P_K}
Let $s>0, \sigma\in\R$, $\alpha\in(-\frac12,2s)$. Then
$\mathcal P_{\!s}$ is a continuous transform
$\mathcal P_{\!s}:\mathcal H^\sigma_\mathcal L\to H^{\alpha+\frac12}(\R\to \mathcal H^{\sigma-\alpha}_\mathcal L)$.
Moreover,
\begin{equation}
\label{eq:iso2}
\llbracket \mathcal P_{\!s}[u]\rrbracket_{H^{\alpha+\frac12}(\R\to \mathcal H^{\sigma-\alpha}_\mathcal L)}^2= 
\frac{\Gamma(s+\frac12 )^2}{\Gamma(s)^2}~\!\frac{ \Gamma(\alpha+1)   \Gamma(2s-\alpha)}{s\Gamma(2s)}~\!\|u\|_{\mathcal H^\sigma_\mathcal L}^2~\!.
\end{equation}
\end{Theorem}

\proof
Thanks to (\ref{eq:iso1}) we already know that 
$$
\|\mathcal P_{\!s}[u]\|^2_{L^{2}(\R\to \mathcal H^{\sigma-\alpha}_\mathcal L)}
\le \|\psi_s\|^2_{L^{2}(\R)}~\!\|u\|^2_{\mathcal H^{\sigma-\alpha-\frac12}_{\mathcal L}}
\le \lambda_1^{-\alpha-\frac12} \|\psi_s\|^2_{L^{2}(\R)}~\!\|u\|^2_{\mathcal H^{\sigma}_{\mathcal L}}
$$
for any $u\in \mathcal H^\sigma_\mathcal L$, which gives the continuity of 
$\mathcal P_{\!s}:\mathcal H^\sigma_\mathcal L\to L^2(\R\to \mathcal H^{\sigma-\alpha}_\mathcal L)$, as $\lambda_1>0$.

Next, take $u=\sum_j u_j\f_j\in \mathcal H^\sigma_\mathcal L$. By the rescaling properties of the Fourier 
transform we have 
$$
\reallywidehat{\mathcal P_{\!s}[u]}(\xi)=
\sum_{j=1}^\infty  \widehat{\psi_{s,\lambda_j}}(\xi)~\!u_j\f_j=
\sum_{j=1}^\infty \lambda_j^{-\frac12} \widehat{\psi_{s}}(\lambda_j^{-\frac12}\xi)~\!u_j\f_j.
$$
This readily gives

$$
\begin{aligned}
\llbracket \mathcal P_{\!s}[u]\rrbracket_{H^{\alpha+\frac12}(\R\to \mathcal H^{\sigma-\alpha})}^2&=
\sum_{j=1}^\infty \lambda_j^{{\sigma-\alpha}-1}u_j^2\intr|\xi|^{2\alpha+1}\big| \widehat{\psi_{s}}(\lambda_j^{-\frac12}\xi)\big|^2~\!d\xi
\\
&=
\big(\intr|\xi|^{2\alpha+1}\big| \widehat{\psi_{s}}(\xi)\big|^2~\!d\xi\big)\sum_{j=1}^\infty \lambda_j^{\sigma}u_j^2,
\end{aligned}
$$
which proves (\ref{eq:iso2}). This concludes the proof by  Proposition \ref{P:psi_transform}
and Lemma \ref{L:P_smooth}.
\QED

We conclude by stating the next  immediate consequence of Theorems \ref{P:jump} and \ref{T:P_K}, which is related to some 
results in \cite{2015}.

\begin{Corollary}
Let $s>0$. For any $u\in \mathcal H^s_\mathcal L$ it holds that
$$
\begin{aligned}
\|\mathcal P_{\!s}[u]\|_{L^2(\R\to \mathcal H^{s+\frac12}_\mathcal L)}^2&= 
\frac{\sqrt\pi \Gamma\big(2s+\frac12 \big)\Gamma\big(s+\frac12 \big)^2}
{s\Gamma(2s)\Gamma(s)^2}~
 \!\|u\|_{\mathcal H^s_\mathcal L}^2\\
\llbracket \mathcal P_{\!s}[u]\rrbracket_{H^{s+\frac12}(\R\to \mathcal H)}^2&= 
\frac{\Gamma(s+\frac12)^2}{\Gamma(2s)}~\!\|u\|_{\mathcal H^s_\mathcal L}^2.
 \end{aligned}
$$
\end{Corollary}

\end{document}